\documentclass[10pt]{amsart}
\usepackage[utf8]{inputenc}
\usepackage{amsmath,amssymb,amsthm}
\usepackage{enumitem}
\usepackage[hidelinks]{hyperref}
\usepackage{fullpage}
\usepackage{cleveref}
\usepackage{tikz}
\usetikzlibrary{patterns}

\DeclareMathOperator{\slog}{\mathrm{s-log}}

\DeclareMathOperator{\Proj}{\P}
\DeclareMathOperator{\SL}{\mathrm{SL}}
\DeclareMathOperator{\GL}{\mathrm{GL}}
\DeclareMathOperator{\KK}{\mathbb{K}}
\DeclareMathOperator*{\esssup}{ess\,sup}

\newcommand{\N}{\mathbb{N}}
\newcommand{\Z}{\mathbb{Z}}
\newcommand{\R}{\mathbb{R}}
\newcommand{\C}{\mathbb{C}}
\newcommand{\Q}{\mathbb{Q}}

\renewcommand{\P}{\mathbb{P}}
\newcommand{\E}{\mathbb{E}}
\newcommand{\HH}{\mathbb{H}}

\newcommand{\vp}{\varphi}
\newcommand{\ve}{\varepsilon}

\newtheorem{thm}{Theorem}[section]
\newtheorem{fact}[thm]{Fact}

\newtheorem{lem}[thm]{Lemma}
\newtheorem{cor}[thm]{Corollary}
\newtheorem{prop}[thm]{Proposition}

\theoremstyle{definition}
\newtheorem{notation}[thm]{Notation}
\newtheorem{deffo}[thm]{Definition}
\newtheorem{remm}[thm]{Remark}

\numberwithin{equation}{section}

\title[Uniform estimates for random matrix products]{Uniform estimates for random matrix products and applications}

\author{Omar Hurtado and Sidhanth Raman}
\newcommand{\Addresses}{{
  \bigskip
  \footnotesize

  \textsc{Department of Mathematics, University of California, Irvine}\par\nopagebreak
  \textit{E-mail addresses}:  \texttt{ohurtad1@uci.edu}, \texttt{svraman@uci.edu} 
}}

\begin{document}
	
\begin{abstract}
	For certain natural families of topologies, we study continuity and stability of statistical properties of random walks on linear groups over local fields. We extend large deviation results known in the Archimedean case to non-Archimedean local fields and also demonstrate certain large deviation estimates for heavy tailed distributions unknown even in the Archimedean case. A key technical result, which may be of independent interest, establishes lower semi-continuity for the gap between the first and second Lyapunov exponents. As applications, we are able to obtain a key technical step towards a localization proof for heavy tailed Anderson models (the full proof appearing in a companion article), and show continuity/stability (taking the geometric data as input) of various statistical data associated to hyperbolic surfaces.
\end{abstract}
\maketitle
\setcounter{tocdepth}{1}
\tableofcontents

\section{Introduction}
\subsection{Main abstract results}

Let $\mathbb{K}$ denote a local field, i.e. $\R$, $\C$, or a finite extension of $\Q_p$ or $\mathbb{F}_p(\!(t)\!)$. For $d \in \N$, we let $\GL_d(\KK)$ and $\SL_d(\KK)$ denote the general linear group of dimension $d$ and special linear group of dimension $d$ respectively. Letting $\mu$ be a (Borel) probability measure on $\GL_d(\KK)$, and presuming (as we will throughout the paper) that we have the mild moment condition
\begin{equation}\label{logmoment}
    \int \log(\max\{\|M\|,\|M^{-1}\|\})\,d\mu(M) < \infty,
\end{equation}
it was conjectured by Bellman in \cite{Bellman} that generally the random product of i.i.d. matrices with distribution $\mu$ should behave very much like a sum of i.i.d. commutative variables; this has been borne out to an astounding degree, which we will briefly summarize in Section \ref{Background}. One of the first results in this direction was a law of large numbers, due to Furstenberg and Kesten:
\begin{thm}[\cite{FurstenbergKesten}]
For any Borel probability measure $\mu$ on $\GL_d(\R)$ satisfying (\ref{logmoment}) and i.i.d. matrices $M_i$ with law $\mu$,
\[ \frac{1}{n}\log\|M_n M_{n-1}\cdots M_1\| \rightarrow \lambda_1(\mu)\]
almost surely, where
\[\lambda_1(\mu) := \lim_{n\rightarrow \infty} \frac{1}{n} \E[\log\|M_n M_{n-1} \cdots M_1\|].\]
\end{thm}

Under some natural assumptions on the support of $\mu$, the same law of large numbers (with the same ``mean'' $\lambda_1(\mu)$) also holds for the quantities $\log\|M_n\cdots M_1x\|$ and $\log |f(M_n\cdots M_1x)|$, where $x$ is a non-zero element of $\KK^d$ and $f$ is a non-zero element of the dual $(\KK^d)^\ast$. Moreover, in the breakthrough work of Benoist and Quint in \cite{bqclt}, it was demonstrated that if one has the stronger moment condition

\begin{equation}\label{finitevar}
    \int \log(\max\{\|M\|,\|M^{-1}\|\})^2\,d\mu(M) <\infty
\end{equation}
and some natural assumptions on the support of $\mu$ (see Definition \ref{nddcond}) then in fact there is a central limit theorem. 
\begin{notation}
    We use $A_n^\mu$ to mean a variable distributed as $M_n\cdots M_1$, with $M_i$ i.i.d. with law $\mu$.
\end{notation}

\begin{thm}[\cite{bqclt}]\label{clt}
    Let $\mu$ be a distribution satisfying (\ref{finitevar}) and which is strongly irreducible and contracting (see Definition \ref{nddcond}). Then there is a constant $\sigma(\mu) \geq 0$ such that
    \begin{equation}
        \frac{\log\|A_n^\mu\| - n\lambda_1(\mu)}{\sqrt{n}} \rightarrow \mathcal{N}(0,\sigma(\mu)^2)
    \end{equation} where $\mathcal{N}(0,\sigma(\mu)^2)$ is the normal distribution with variance $\sigma(\mu)^2$ and mean zero, and the convergence is in distribution. If $\KK = \R$ or $\C$, then $\sigma(\mu) > 0$.
\end{thm}

Such a theorem had was first proven by Le Page under a stronger moment assumption \cite{LePage}, and later generalized by Jan \cite{janthesis}. (An earlier version of this result was proven by \cite{tutubalin}, but under the stringent assumption of a density with respect to the Haar measure, and for semigroups of positive matrices, Furstenberg and Kesten obtained an earlier central limit theorem.)

Under stronger moment assumptions than (\ref{logmoment}) and the same assumptions on the support of $\mu$, one can obtain large deviation estimates, i.e. estimates on the rate of convergence in the law of large numbers:
\begin{equation}\label{llndecay} \P\left[\left|\frac{1}{n}\log\|M_n\cdots M_1\| - \lambda_1(\mu)\right| > \ve \right] \rightarrow 0\end{equation}
For example, if $\int e^{\alpha \log \|M\|} \,d\mu(M)$ is finite for some $\alpha > 0$, then one has exponential decay in (\ref{llndecay}); this is originally due to Le Page and now classical. If $\int e^{ |\log \|M\| |^\delta}\,d\mu(M) < \infty$ for some $\delta$, one has ``semi-exponential'' decay in (\ref{llndecay}), i.e. an upper bound of the form $Ce^{-cn^\delta}$ for $C,c$ positive constants; this was obtained recently by Cuny, Dedecker, and Merleverde \cite{cuny2017large}. Prior to the work of Cuny et al., Benoist and Quint showed that if $\log\max\{\|M\|,\|M^{-1}\|\}$ has a polynomial moment of order $p\geq 1$, then one gets e.g. $o(n^{2-p})$ decay in (\ref{llndecay}) (\cite{cuny2017large} in fact relies on \cite{bqclt}).\footnote{They actually proved the stronger statement that $\sum_n \P\left[\left|\frac{1}{n}\log\|M_n\cdots M_1\| - \lambda_1(\mu)\right| > \ve \right] n^{p-2} < \infty$.}

Our work for the most part concerns the stability of these estimates under small perturbations of the measure $\mu$. The appropriate sense in which to understand what constitutes a ``small perturbation'' is in terms of certain topologies: those induced by Wasserstein metrics and certain generalizations thereof. While we will explicitly define these metrics in Section \ref{prelimsec}, we introduce the important metrics and the spaces on which they are defined somewhat informally here.

When we work with distributions on $\GL_d(\mathbb{K})$, we work with the distance 
$$d(M,M') = \max\{\|M-M'\|,\|M^{-1}-(M')^{-1}\|\};$$ on $\SL_d(\mathbb{\KK})$ it suffices to work with the metric $d(M,M') = \|M-M'\|$. For most of our results, we will treat the $\SL_d(\KK)$ case explicitly; the $\GL_d(\KK)$ makes computations slightly more unwieldy but does not otherwise introduce any new complications. We emphasize that this is only true because of the choice of metric on $\GL_d(\KK)$; if one instead equips $\GL_d(\KK)$ with the metric $d(M,M') = \|M-M'\|$, our results do not hold in general. It is necessary to control both the largest singular value and the inverse of the least singular value; on the special linear group these are related by the fact that they multiply to unity, whereas in the general linear case the finer metric accounts for the least singular value. See \Cref{glvssl} for more details.

\begin{deffo} We now define spaces of probability measures satisfying certain moment conditions:
    \begin{itemize}
        \item We let $\mathcal{P}(\SL_d(\KK))$ denote the space of Borel probability measures on $\SL_d(\KK)$, and similarly for the space of probability measures on $\GL_d(\KK)$.
        \item For $p \geq 1$, let $\mathcal{P}_{\log}^p(\SL_d(\KK))$ denote the subset of $\mathcal{P}(\SL_d(\KK))$ which consists of distributions with logarithmic moments of order $p$, i.e. those for which
        \[ \int_{\SL_d(\KK)} (\max\{0, \log\|M\|\})^p \,d\mu(M) < \infty\]
        \item For $\delta \in (0,1)$, we let $\mathcal{P}_{\slog}^\delta(\SL_d(\KK))$ denote the subset of $\mathcal{P}(\SL_d(\KK))$ which has a semi-logarithmic moment of order $\delta$, i.e. those for which 
        \[ \int_{\SL_d(\KK)} \exp(\max\{\log \|M\|, 0 \}^\delta) \,d\mu(M) < \infty\]
        \item For $\alpha > 0$, we let $\mathcal{P}^\alpha(\SL_d(\KK))$ denote distributions with a moment of order $\alpha$, i.e. those for which
        \[ \int_{\SL_d(\KK)} \|M\|^\alpha \,d\mu(M) < \infty\]
        \item For all of these, we define $\mathcal{P}^p_{\log}(\GL_d(\KK))$ and so on analogously, except that in the general linear case one must replace $\|M\|$ with $\max\{\|M\|,\|M^{-1}\|\}$.
    \end{itemize}
\end{deffo}

We will make use of Wasserstein (or Wasserstein-like) metrics which are quite natural for these spaces. As mentioned before, we will postpone the precise definitions to \Cref{prelimsec}, but the following facts are straightforward consequences of existing work on Wasserstein metrics.
\begin{prop}
    The families of metrics $W_{\log}^p$, $W_{\slog}^\delta$, and $W^\alpha$, defined on the spaces $\mathcal{P}_{\log}^p(\SL_d(\KK))$, $\mathcal{P}_{\slog}^\delta(\SL_d(\KK))$, and $\mathcal{P}^\alpha(\SL_d(\KK))$ respectively, induce topologies characterized by the following:
    \begin{itemize}
        \item $W^p_{\log}(\mu_n,\mu) \rightarrow 0$ if and only if $\mu_n \rightarrow \mu$ in the weak topology and moreover the $p$-th logarithmic moments of $\mu_n$ converge to that of $\mu$.
        \item $W^\delta_{\slog}(\mu_n,\mu) \rightarrow 0$ if and only if $\mu_n \rightarrow \mu$ weakly and moreover the $\delta$ order semi-logarithmic moments of $\mu_n$ converge to that of $\mu$.
        \item $W^\alpha(\mu_n,\mu) \rightarrow 0$ if and only if $\mu_n \rightarrow \mu$ weakly and morevoer the $\alpha$-th order moments of $\mu_n$ converge to those of $\mu$.
    \end{itemize}
\end{prop}

For $\alpha \geq 1$, the distances $W^\alpha$ have been studied quite widely for other metric spaces, for e.g. the study of optimal transport (see e.g.\cite{villani2009optimal} and \cite{figalli2021invitation}) and in computer vision (see e.g. \cite{Rubner2001}). The $W^\alpha$ distances for $\alpha < 1$ have been studied already in the context of random matrix products by \cite{tall2020moduli}. Our main abstract results are essentially ``stability'' of various existing estimates; for example the first set of results essentially say that if one perturbs a measure a small bit one does not have blow up of the large deviation constants. Consequently, and more useful for applications, we can extract uniform estimates over compact sets of distributions in the appropriate topologies.

\begin{thm}\label{logmatrixlde}
	Let $K\subset \mathcal{P}^{p}_{\log}(\SL_d(\mathbb{K}))$ be a set of probability measures such that
	\begin{enumerate}
		\item $K$ is compact in the $W^{p}_{\log}$ topology,
		\item Every $\mu \in K$ is strongly irreducible and contracting (see Definition \ref{nddcond}).
	\end{enumerate}
	Then there are constants $C_n = C_n(K,\ve)$ such that $\P[ |\log\|A_n^\mu x\| - n\lambda_1(\mu)|>n\ve]\leq C_n$ and moreover, these $C_n$ satisfy
    \begin{equation}\label{sumcondition}
	\sum_{n\in\N} n^{p-2} C_n < \infty
    \end{equation}
\end{thm}
\begin{remm}
    One can improve on (\ref{sumcondition}) and obtain the following:
    \begin{equation}
        C_n \leq \begin{cases} Cn^{\frac{3}{2} - \frac{3}{2}p} \quad &\text{for} \quad 1 < p \leq 3\\
        Cn^{-p}\quad &\text{for} \quad p > 3
        \end{cases}
    \end{equation}
    for some fixed $C$ depending on $K$ and $\ve$. These stronger bounds are obtained by carefully tracing through the estimates that appear in the proof of \Cref{explicitmgpolylde}. We note that such bounds are essentially due to Benoist and Quint, and record them here to use them when optimizing \Cref{semiloglocal}. We wrote \Cref{logmatrixlde} in terms of (\ref{sumcondition}) to emphasize that it is a uniform version of \cite[Proposition 4.1]{bqclt}.
\end{remm}

\begin{thm}\label{semilogmatrixlde}
	Let $K\subset \mathcal{P}^{\delta}_{\slog}(\SL_d(\mathbb{K}))$ be a set of probability measures such that
	\begin{enumerate}
		\item $K$ is compact in the $W^{\delta}_{\slog}$ topology,
		\item Every $\mu \in K$ is strongly irreducible and contracting.
	\end{enumerate}
	Then there are constants $C=C(K,\ve,\delta)$ and $c=c(K,\ve,\delta)$ such that
		\[\P[|\log\|A_n^\mu x\| - n\lambda_1(\mu)|>n\ve]\leq Ce^{-cn^{.99\delta}} \]
	for any $x \in \mathbb{K}^d$ with $\|x\|=1$. 
\end{thm}
\begin{remm}
    Here and throughout, $.99\delta$ can be replaced by any $\delta' < \delta$, with the constants appearing then gaining a dependence on $\delta'$. To reduce notational load, we just write $.99\delta$.
\end{remm}
We also produce large deviation estimates for matrix coefficients, though our argument requires technical assumptions on the image of $K$ under the pushforward of the second exterior power in the weakest moment regime:
\begin{thm}\label{matrixcoefflde}
In the settings of Theorems \ref{logmatrixlde} and \ref{semilogmatrixlde}, so long as $\wedge^2_{\ast}(K)$ is compact in  $\mathcal{P}^1_{\log}(\SL_{\binom{d}{2}}(\mathbb{K}))$, then said theorems remain true if one replaces $\log\|A_n^\mu x\|$ with $\log |f(A^\mu_nx)|$, for any $x\in \mathbb{K}^d$ and $f\in (\mathbb{K}^d)^\ast$ with $\|x\|=\|f\|=1$. This compactness assumption is automatically satisfied if $K$ is compact in $\mathcal{P}^\delta_{\slog}$ for some $\delta \in (0,1)$.
\end{thm}

We can also prove exponential large deviations on compact sets in $\mathcal{P}^\alpha(\mathrm{SL}_d(\mathbb{K}))$ for $\mu$ which satisfy the dynamical conditions and possess some fractional moment:

\begin{thm}\label{unifldeexp}
    Let $K \subset \mathcal{P}^\alpha(\mathrm{SL}_d(\mathbb{K}))$ be a set of probability measures such that
    \begin{enumerate}
        \item $K$ is compact in the $W^\alpha$ topology,
        \item Every $\mu \in K$ is strongly irreducible and contracting.
    \end{enumerate}
    Then there are constants $C= C(K,\ve)$ and $c= c(K,\ve)$ such that
    \begin{equation}\label{expldeseqn} \P\left[|\log \|A_n^\mu x\| - n\lambda_1(\mu)|>n\ve \right] \leq Ce^{-cn} \end{equation}
    for any $x \in \mathbb{K}^d$ with $\|x\| = 1$.
\end{thm}

As we will discuss at further length shortly, there are results known which are considerably stronger than ours in special cases \cite{duarte2016lyapunov, duarte2020large}. To the best of our knowledge, this result is new at this level of generality, even for Archimedean fields $\mathbb{K} = \R$ or $\mathbb{C}$, and we are not aware of any results of this type for non-Archimedean fields.

We quickly remark that there is an intrinsic interest in the study of random walks on semisimple groups over non-Archimedean local fields from more algebraic and number theoretic disciplines. For example, the study of random products of matrices over $\Q_p$ has led to central limit theorems for random walks on Bruhat--Tits buildings \cite{CartwrightWoess,parkinson2007isotropic}; we do not specifically pursue the consequences of our stability results in this context.

In the polynomial moment case, the uniformity for large deviations of the matrix elements is novel. We believe large deviation estimates for the matrix elements in the semi-exponential moment regime to be novel, even before one accounts for the uniformity; however, for the quantities $\|A_n^\mu x\|$ and $\|A_n^\mu\|$ the non-uniform version of our estimates was obtained in \cite{cuny2017large}.

We also demonstrate the continuity of various statistical data in these topologies when restricting to the locus of measures which are strongly irreducible and contracting. The first result is essentially an immediate consequence of work of Furstenberg and Kifer in \cite{FurstenbergKifer} together with Chebyshev; we include it to emphasize that the logarithmic Wasserstein metrics are in fact quite natural.

\begin{prop}\label{lyapcont}
    The function $\mu \mapsto \lambda_1(\mu)$ is continuous in the $W^1_{\log}$ topology when restricted to the locus of measures which are strongly irreducible and contracting, and necessarily in all the finer Wasserstein-type topologies.
\end{prop}

We also obtain continuity of the variance appearing in the central limit theorem of Benoist and Quint, \cite[Theorem 1.1]{bqclt}. This was non-trivial, and required establishing continuity in the datum $\mu$ of a certain coboundary term constructed by Benoist and Quint in \cite{bqclt} to prove their central limit theorem.

\begin{thm}\label{varcont}
    The function $\mu \mapsto \sigma(\mu)$ is continuous in the $W^\delta_{\slog}$ topology, where $\sigma(\mu)$ is the standard deviation appearing in \Cref{clt}.
\end{thm}

To our knowledge, continuity of the variance is completely new; we note that in analogy with the scalar case, one would expect continuity in the weaker $W^2_{\log}$ topologies; we obtain the analogue of our result in this weaker topology in the special case of $\SL_2(\KK)$ cocycles:

\begin{thm}\label{varcont2d}
    The function $\mu \mapsto \sigma(\mu)$ is continuous as a map $W^2_{\log}(\SL_2(\KK)) \rightarrow \R_{\geq 0}$ when one restricts to any locally compact subset of the locus of strongly irreducible and contracting measures.
\end{thm}

 The difference between the case $\SL_2(\KK)$ and other cases is related to a technical result on stability of the ``gap'' $\lambda_1(\mu) - \lambda_2(\mu)$, our \Cref{wedgecontinuity}; see also \Cref{wedgecontrmk}. It is a general theme in the study of linear cocycles over various base dynamics that the leap from $d=2$ to $d>2$ is quite difficult.

Though it is somewhat straightforward given our large deviation estimates, we mention that  combining many of the estimates used to obtain these results with the argument used to prove a regularity result by Benoist and Quint in  \cite[Proposition 4.5]{bqclt}, we ourselves obtain various intermediate regularity result for invariant measures coming from $\mu$ with semi-logarithmic moments. (We recall that any $\mu$ on $\SL_d(\KK)$ which is strongly irreducible and contracting has a unique invariant measure $\nu$ on $\P(\KK^d)$.) We present one result of this flavor here to give an idea of what results we obtain. This requires the introduction of a few notions; for the following we assume $\R^d$ and $\C^d$ have the usual inner product structure, and define a distance $d$ on the projectivizations $\P(\R^d)$, $\P(\C^d)$ by:
\[ d(x,y)^2 = 1 - |\langle x, y \rangle|^2\]
where $x,y$ abusively denote both equivalence classes of vectors under non-zero scaling and specific normalized representatives. Note that $d(x,y)$ is also the sine of the minimal angle between two representatives of the equivalence classes.

We recall a notion used in e.g. \cite{duarte2020large}:
\begin{deffo}
    We call a measure $\nu$ on $\P(\KK^d)$ (for $\KK = \R$ or $\C$) \emph{weak-H\"older continuous} if there are $c,\rho > 0$ such that sufficiently small $\ve > 0$ and any $x \in \P(\KK^d)$ we have \[\nu(B_\ve(x)) \leq \exp(-c(\log(\ve^{-1})^\rho)\]
    where $B_\ve(x)$ is the open ball of radius $\ve$.
\end{deffo}
\begin{remm}
     The notion of weak-H\"older in \cite{duarte2020large} is defined for functions; much like the notion of H\"older continuity for measures is a straightforward generalization of the notion for functions, this definition is the analogue for measures of the notion therein. We also mention that while the name weak-H\"older was not used, this intermediate modulus of continuity which ``interpolates'' between H\"older continuity and log-H\"older continuity (a notion we will discuss later) has appeared in the study of linear cocycles at least as far back as the work of Goldstein and Schlag on regularity of the Lyapunov exponents for certain quasi-periodic cocycles \cite{goldstein2001holder}.
\end{remm}
A representative regularity result is the following:
\begin{thm}\label{regularitythmintro}
Given any measure $\mu \in W^\delta_{\slog}(\SL_d(\KK))$ for $\KK= \C$ or $\R$ which is moreover strongly irreducible and contracting, the associated invariant measure $\nu$ is weak-H\"older continuous.
\end{thm}
This is a consequence of a more general regularity result, \Cref{semilogfunctionalreg}.
It is now classical, originally due to Le Page, that if one has exponential moments, one obtains the stronger property of H\"older continuity of the invariant measure, and in a certain sense log-H\"older continuity (we will define and discuss this in \Cref{invsec}) of the invariant measure is the essential technical result in \cite{bqclt}.

Finally, a key technical result which may be of independent interest concerns the gap between the first and second Lyapunov exponents. The first Lyapunov exponent $\lambda_1$ is essentially the rate of exponential growth for the first singular value of random matrix products. The second Lyapunov exponent $\lambda_2$ is essentially the rate of exponential growth of the second singular value (see \Cref{rmp} for a precise definition).  

\begin{thm}\label{unifgapbound}
    Let $K$ be a compact subset of $\mathcal{P}_{\slog}^\delta$ such that every $\mu \in K$ is strongly irreducible and contracting. Then
    \begin{equation}
        \inf_{\mu \in K} \lambda_1(\mu) - \lambda_2(\mu) > 0
    \end{equation}
\end{thm}

Pointwise positivity for $\mu$ strongly irreducible and contracting was shown by Guivarc'h; this result essentially follows from \Cref{lyapcont} together with the following result:

\begin{thm}\label{uscsum}
    The map $\mu \mapsto \lambda_1(\mu) + \lambda_2(\mu)$ is upper semi-continuous.
\end{thm}

While such results have been obtained in certain contexts for either one parameter families or compactly supported measures $\mu$, this seems to be the first result of its kind for unbounded distributions at this level of generality; we discuss existing results in more detail in the next section.

\subsection{Background on random matrix products}\label{Background}

The study of random matrix products is a mature subject at this point; for a comprehensive introduction we recommend the books \cite{bougerol2012products, benoist2016random}. We will nevertheless sketch a broad account of the history, inevitably biased towards those results which clarify the ways ours fit into the broader theory.

The pointwise versions of the theorems we are concerned with, at least in the regime where we have exponential moments for $\log\|M\|$ (or equivalently, power moments for $\|M\|$), are now classical.\footnote{We also presume strong irreducibility and contraction.} The law of large numbers was established by Furstenberg and Kesten, and the Central Limit Theorem and Large Deviations both by Le Page. The same work of Furstenberg and Kesten in fact proved a CLT for semigroups of positive matrices, and the work of Le Page had some additional technical hypotheses shown to be unnecessary in \cite{guivarc1986products,gol1989lyapunov}. The work of Furstenberg and Kesten was in fact optimal, at least as far as moment conditions are concerned.

Absent the assumption of compact support, the question of continuity of the top Lyapunov exponent seems to have been more or less settled by Furstenberg and Kifer, which among other things demonstrated the following:
\begin{thm}[\cite{FurstenbergKifer}]
    Let $\mu$ be a distribution on $\GL_d(\R)$ which is strongly irreducible and contracting, satisfying
    \[ \int \log \max\{\|M\|,\|M^{-1}\|,0\}\,d\mu(M) < \infty.\]
    If $\mu_n\rightarrow \mu$ weakly and satisfies
    \[ \sup_{n\in\N} \int_{\|M\| > T} \log\max\{\|M\|,0\}\,d\mu_n(M) + \int_{\|M^{-1}\| > T} \log\max\{\|M^{-1}\|,0\}\,d\mu_n(M) \rightarrow 0,\]
    as $T \to \infty$, i.e. the tails decay uniformly, then we have convergence of the top Lyapunov exponents $\lambda_1(\mu_n) \rightarrow \lambda_1(\mu)$. 
\end{thm}
In fact, this is slightly weaker than what they proved; strong irreducibility and contraction are sufficient conditions for the various phenomena we are concerned with throughout the paper, but not always necessary. We are concerned with distributions which are not necessarily compactly supported, where the continuity question has largely not progressed since the work of Furstenberg and Kifer, at least not qualitatively; for certain one parameter families with exponential moment conditions, H\"older regularity was established by Le Page for the top Lyapunov exponent \cite{page1989regularite}. We note that within the space of compactly supported measures equipped with an appropriate topology, there has been significant investigation recently, qualitative and quantitative, regarding the question of continuity --- see e.g. \cite{avila2023continuitylyapunovexponentsrandom,tall2020moduli,malheiro2015lyapunov,bocker2017continuity,duarte2016lyapunov,duarte2020large}.

While continuity of the Lyapunov exponent is something which is very well understood in the setting where $\mu$ is strongly irreducible and contracting, there is recent progress on quantitative questions related to these phenomena. For example, there has been work getting precise asymptotics for the large deviation estimates for distributions with exponential moments in \cite{xiao2020precise}. Recently Berry--Esseen estimates were obtained, first for the exponential moment regime \cite{xiao2022berry} and then for polynomial moments \cite{cuny2023berry}. In a pointwise sense, corresponding large deviations were obtained essentially by \cite{bqclt}, and then in the semi-exponential moment regime by \cite{cuny2017large}.

Uniform results are generally constrained to the bounded (i.e. compactly supported) case, with the exception of \cite{tsay1999some} which proved uniform estimates for one parameter families (see also \cite{Bucaj2017LocalizationFT} for another proof, though only for the bounded case). To our knowledge, all uniform large deviation type results in the literature, including this work of Tsay, are constrained to the Archimedean case. In \cite{duarte2016lyapunov}, Duarte and Klein use spectral methods inspired by earlier work of Le Page \cite{LePage} and Hennion and Herve \cite{hennion2001limit} to obtain uniform large deviations for cocycles understood as deterministic maps over sufficiently random dynamical systems, those called strongly mixing Markovian (see \cite[Theorem 5.2]{duarte2016lyapunov} for a precise statement). In particular, special cases of their work include e.g. Schr\"odinger cocycles corresponding to a bounded random potential and finitely supported distributions corresponding to a fixed probability vector. Both of these correspond to Bernoulli shifts equipped with an ergodic measure (the former with possibly infinite alphabet). The work in \cite{duarte2016lyapunov} required irreducibility hypotheses, but in \cite{duarte2020large} the same authors were able to get uniform (subexponential) large deviation estimates for $\mu$ which were not contracting or strongly irreducible in the special case of finitely supported distributions corresponding to a fixed probability vector. As a consequence, the authors were able to obtain local modulus of continuity results for the Lyapunov exponents using a general approach introduced in \cite{duarte2016lyapunov}. 

Besides the ability to treat more general Markovian systems, the most striking advantage of the work of Duarte and Klein over ours is that they can precisely pin down the dependence on $\ve$ of the constants appearing in \Cref{expldeseqn}, e.g. $c$ scales as $\ve^2$, in analogy with the scalar case, where Hoeffding's inequality yields precisely this scaling for large deviations of sums of bounded random variables. At the present time, our methods do not yield any information about the asymptotic dependence of constants on $\ve$, due to the use of compactness arguments, though we expect that at least in the subexponential regime it should be possible to replace these compactness arguments with quantitative estimates.

A method introduced in \cite{Bucaj2017LocalizationFT} was specifically meant to treat Schr\"odinger cocycles, but appears to be fairly robust as long as one works with compactly supported cocycles, which are strongly irreducible and contracting. Besides the aforementioned work of Tsay, all of these results require almost sure boundedness.

We briefly remark that essentially all the work here concerns i.i.d. random matrix products; a family of questions very similar in spirit (though quite different in technical details) is that of non-stationary random matrix products. This was investigated by Gorodetski and Kleptsyn \cite{gorodetski2022non}, who were able to prove results similar to the i.i.d. case under reasonable assumptions, i.e. the existence of something like a Lyapunov exponent and large deviation estimates. Very recently, these authors and Monakov proved a non-stationary central limit theorem in \cite{gorodetski2024central}.

Finally, regarding uniform estimates on $\lambda_1(\mu) - \lambda_2(\mu)$, we mention first that in the setting of compactly supported measures under the topology used in e.g. \cite{bocker2017continuity,avila2023continuitylyapunovexponentsrandom}, an analogue of \Cref{unifgapbound} is very straightforward; the strategy is essentially that carried out here, except that the analogue of a key technical result (\Cref{wedgecontinuity}) is much easier to prove for compactly supported distributions. For various one-parameter families of interest in mathematical physics, one can similarly obtain ``uniform simplicity'' of the Lyapunov spectrum. This is crucial in certain work on generalizations of the Anderson model to ``quasi one-dimensional'' contexts e.g. \cite{klein1990localization,Macera2022}. However, such results require quite strong geometric assumptions on $\mu$ which go beyond just strong irreducibility and contraction (see also \cite{gol1989lyapunov}).

\subsection{Random Schr\"odinger operators}
The theory of random matrix products has had a longstanding and fruitful intersection with the study of random Schr\"odinger operators in one dimension; via the transfer matrix formalism, one can study the asymptotics of eigenvectors of such operators via such random products.

Specifically, we are concerned with operators acting on $\ell^2(\Z)$ of the form
\begin{equation}\label{tischrodeq}
    H = \Delta + V
\end{equation}
where $\Delta$ is a discrete Laplacian, defined by $[\Delta \psi](n) = \psi(n+1)+\psi(n-1)$ for $\psi \in \ell^2(\Z)$, and $V$ is a random potential acting by multiplication, $[V\psi](n) = V_n\psi(n)$, with $V_n$ taken independent and identically distributed with law $\mu$. (It is more physically natural to replace $\Delta$ here with the negative of the graph Laplacian when $\Z$ is treated as a graph with integers distance exactly one apart connected; the spectral theory is exactly the same and our chosen formulation simplifies many computations.)

One of the central questions in the study of random Schr\"odinger operators is localization; this is, in physical terms, a phenomenon whereby disorder (the randomness) thwarts transport, and electrons essentially become trapped. The strongest forms of localization correspond to strong bounds on the moments of the position operator composed with the time evolution operator $e^{-itH}$, and various notions like these are variously called dynamical localization. A variant which is in general strictly weaker (see \cite{del1996operators}) but is often a strong indicator of dynamical localization is the Anderson localization, which corresponds to the absence of a continuous spectrum and exponential decay of the eigenvectors of $H$. Anderson localization, sometimes called exponential localization, does imply some bounds on the quantum dynamics i.e. on moments of the position operator as the system evolves in time, but these are weak.

We formulate Anderson localization more precisely:
\begin{deffo}
    An operator $H$ is \emph{Anderson localized} if it does not have any continuous spectrum and moreover its associated eigenfunctions decay exponentially. We say $H$ is Anderson localized in $I \subset \R$ if it has no continuous spectrum in $I$ and the associated eigenfunctions with energies $E \in I$ decay exponentially.
\end{deffo}

We will not try to summarize this field of study in its entirety, as it is quite mature; we recommend e.g. \cite{aizenman2015random,cycon2009schrodinger} for broad accounts of the theory. While there has been significant progress in dimensions higher than one, for regular potentials since the eighties (see e.g. \cite{frohlich1985constructive,aizenman1993localization,von1989new,damanik2001multi}) and for singular potentials since the aughts (see e.g. \cite{bourgain2005localization,klein2012comprehensive,li2022anderson,li2022anderson2,hurtado2024localization}), results on higher dimensional lattices require entirely different methods. In one dimension, the very first methods were in some sense ad hoc and required regularity of the noise in addition to boundedness \cite{goldsheid1977random,kunz1980spectre}. Combining the multi-scale analysis of Fr\"ohlich and Spencer with estimates coming from the transfer matrix approach, Carmona, Klein, and Martinelli were able to prove localization for distributions which were not necessarily regular under a fairly mild moment assumption:
\begin{thm}[{\cite{carmona1987anderson}}]
	Let $\mu$ be a non-trivial probability measure on $\R$, i.e. one supported on at least two points. If there exists $\alpha > 0$ such that $\int |x|^\alpha \,d\mu(x) <\infty$, then the random Schr\"{o}dinger operator $H= \Delta + V$ (where $V_n$ are i.i.d. with law $\mu$) is Anderson localized throughout the spectrum, almost surely.
\end{thm}

Since this first proof, there have been many subsequent proofs which have avoided using the multi-scale analysis techniques (e.g. \cite{Bucaj2017LocalizationFT,Jitomirskaya2019, gorodetski2021parametric, rangamani2019singular,gorodetski2024non}), proofs which could be called ``purely one-dimensional'', relying principally on the asymptotics of random matrix products and resultant asymptotics for generalized eigenvectors. While this is not quite a theorem, it is generally believed that ``uniform large deviations plus positive Lyapunov exponent implies localization''. In particular, there are no known models which exhibit both these properties but not localization.

By large deviations we mean large deviation estimates for the so-called transfer matrices, which encode the asymptotics of formal solutions to the eigenequation. We will discuss the significance of these transfer matrices further in Section \ref{schrodsec}; for any potential $V_n \in \R^{\Z}$ and energy $E \in \R$ the transfer matrices associated to the corresponding Schr\"odinger operator are defined as follows:

\begin{equation}
    A_n^E := \begin{pmatrix} E - V_n & -1 \\ 1 & 0 \end{pmatrix} \begin{pmatrix} E - V_{n-1} & - 1 \\ 1 & 0 \end{pmatrix}\cdots \begin{pmatrix} E - V_1 & -1 \\ 1 & 0 \end{pmatrix} 
\end{equation}

It is well known, by many of the results discussed in the background on random matrix products, that if the $V_n$ are i.i.d. with law $\mu$ satisfying a very mild moment condition, then for all $E \in \R$ there is some $\lambda(E) > 0$,  the Lyapunov exponent for energy $E$,  such that
\[
\frac{1}{n} \log \|A_n^E\| \rightarrow \lambda(E)
\]
almost surely.

Our results give a quantitative version of this result for $\mu$ satisfying the following with $p > 1$:
\begin{equation}\label{semilogbd}
    \int \max\{\log x, 0 \}^p\,d\mu(x) < \infty.
\end{equation}
In particular, the following is a consequence of our uniform estimates:

\begin{thm}\label{semilogschrodldes}
    Let $\mu$ be a measure satisfying the bound (\ref{semilogbd}) for some $p > 1$, which is supported on at least two points. Then for any $I \subset \R$ compact, and $\ve > 0$, there is $C = C(I,\ve, \mu)$ such that
    \begin{equation}
        \P\left[\left|\frac{1}{n}\log \left|\langle y, A_n^E x \rangle\right| - \lambda(E)\right| > \ve\right] \leq C\max\{n^{\frac{3}{2}-\frac{3}{2}p}, n^{-p}\}
    \end{equation}
    for any $E \in I$, $\|x\| = \|y\| = 1$.
\end{thm}

These large deviation estimates are indeed enough to prove localization for sufficiently large $p$, though the ``purely'' one-dimensional methods seem not to suffice. Instead, using the estimates from above, it is possible to verify the hypotheses for a multi-scale analysis, as in e.g. \cite{carmona1987anderson, von1989new}, and obtain localization.

In order to avoid going through the precise details of what the hypotheses are and introducing all the notions necessary to discuss a multi-scale analysis in this paper, the full proof of localization appears in a note by the first named author \cite{hurtado2025localization}, where the following is proven:

\begin{thm}[\cite{hurtado2025localization}]\label{semiloglocal}
    Let $\mu$ be any measure supported on at least two points which satisfies (\ref{semilogbd}) for some $p>11$. Then the operator $H = \Delta + V$ for $V_n$ i.i.d. with law $\mu$ is almost surely Anderson localized. 
\end{thm}

It is worth emphasizing that the original proof of Carmona, Klein, and Martinelli and all other listed work is capable of treating the case of singular distributions, e.g. Bernoulli, for the potential. Under certain regularity conditions, the result was known earlier. Moreover, though this was after the work of Carmona, Klein, and Martinelli, given the assumption of absolute continuity of the measure $\mu$ with respect to Lebesgue measure and a bounded density, the fractional moment method developed by Aizenman and Molchanov obtains the following:
\begin{thm}[\cite{aizenman1993localization}]
    If $\mu$ is of the form $d\mu = f\,dx$ for $f$ a bounded density function, then the operator $H = \Delta + V$ where the $V_n$ are i.i.d. with law given by $\mu$ is almost surely Anderson localized.
\end{thm}

In particular, there is no moment condition, and so the forthcoming localization result is only novel in the regime where $\mu$ is singular in some way, as otherwise the result can be proven by the fractional moment method.

While our large deviation results will be used to obtain a localization result in the forthcoming \cite{hurtado2025localization}, it is worth mentioning briefly that Wasserstein distances have already been used fruitfully in the study of random Schr\"odinger operators, specifically in the study of an associated object called the \emph{density of states measure}. The density of states (DOS) measure is of great importance physically, and is roughly an asymptotic normalized eigenvalue count. For the one-dimensional Anderson model associated to $\mu$, Hislop and Marx \cite{hislop2020dependence1,hislop2020dependence2,hislop2021dependence} and Shamis \cite{shamis2021continuity} have studied the way in which the DOS measure associated to $\mu$ depends on $\mu$, and under various assumptions, demonstrated continuity in $\mu$ with respect to certain Wasserstein distances.

\subsection{Statistics of random geodesics on hyperbolic surfaces}

The study of statistical properties of random processes on Riemannian manifolds have led to a deeper understanding of their function theory and global geometric structures. Specializing to hyperbolic geometry, various problems such as random displacement distance for convex-cocompact actions on $\HH^2$ \cite{pollicott1998comparison}, counting self-intersection of oriented curves on surfaces \cite{chas2012self}, and counting random closed geodesics of bounded length on hyperbolic surfaces \cite{gekhtman2019central} all enjoy central limit theorems. 

While our work is not applicable in the context of the results of Gekhtman, Taylor, and Tiozzo, it allows us to make a first step and provides evidence that the statistical data does in fact depend continuously on the geometric structure of the surface. In particular, the central limit theorems proved by Gekhtman, Taylor and Tiozzo relate the \emph{algebraic length} of a geodesic (the length of a minimal word with respect to a generating set) to its geometric length. While there is a straightforward relationship between geometric length and the norm of certain random matrices, these matrices are not exactly distributed like random matrix products when one samples uniformly among those homotopy classes with algebraic length $n$.

If one instead one samples randomly among all words of length $n$, then the distributions of random matrices one considers can in fact be studied using random matrix products, and one also obtains a central limit theorem; see e.g. \cite{park2019probabilitylawsdistributiongeometric}. For this type of counting, our uniform estimates for random matrix products further implies that the parameters (mean, variance) appearing in the central limit theorems depend continuously on the geometric data in a certain sense.

To state our theorem, let $(\Sigma,g) = \HH^2/\Gamma$ be a finite type hyperbolic surface, where $\Gamma < \mathrm{PSL}_2(\R)$ is given the standard surface group presentation (or the standard free group presentation if $\Sigma$ is not closed). Let $\mathcal{F}_n$ denote the set of conjugacy classes of words in $\Gamma$ of length $n$ and let $\mu_n$ be the uniform distribution on $\mathcal{F}_n$. Define $\ell_{g}(\gamma)$ to be the length of the unique geodesic representative of the word $\gamma\in \Gamma$. It is well known that if the $A\in \mathrm{PSL}_2(\R)$ is the matrix corresponding to the loop $\gamma$, then $\ell_g(\gamma) = \log \| A\|$.

First we recall the following central limit theorem, which is essentially an immediate consequence of the classical theory of random matrix products together with fundamental results in hyperbolic geometry; an explicit proof in more generality is contained in e.g. \cite{park2019probabilitylawsdistributiongeometric}, though this more general fact is also a consequence of work in \cite{benoist2016random}:
\begin{thm}\label{hyperbolicthm}
    Let $(\Sigma,g)$ be a closed hyperbolic surface, and let $\mathcal{S} \subset \pi_1(\Sigma)$ be the standard symmetric generating set. If $\mu_n$ is the uniform measure on $\mathcal{F}_n$ (the set of conjugacy classes of length $n$ nonreduced words in the alphabet $\mathcal{S}$), then there exists constants $L_g, \sigma_g > 0$ such that
    $$\mu_n\left(\gamma : \frac{\ell_g (\gamma)  - n L_g}{\sigma_g\sqrt{n}}\in [a,b]\right) \to \frac{1}{\sqrt{2\pi}}\int_a^b e^{-x^2/2} dx.$$
\end{thm}
Our result is the following:
\begin{thm}\label{HyperbolicCountingCLT}
    Given any closed surface $\Sigma$ admitting some hyperbolic metric, and letting $T(\Sigma)$ be its Teichm\"uller space, i.e. the space of hyperbolic structures, the quantities $L_g$ and $\sigma_g$ are continuous in the topology induced by the Teichm\"uller metric.
    
\end{thm}

We will recall the details of Teichm\"uller space in Section \ref{geodesicsec}. 
Continuous dependence of the central limit on $L_{g}$ easily follows from Furstenberg--Kifer. Continuity of the variance $\sigma_{g}$ appears to be new, and immediately follows from a simple lemma in hyperbolic geometry and our general results on continuity of variance. In fact, we prove an analogous result for finite type hyperbolic surfaces which relies on the fact that parabolic elements in the fundamental group are probabilistically sparse; see \Cref{FiniteTypeHyperbolic}.

Our methods as they stand cannot really be adapted to prove continuity in the results of \cite{gekhtman2019central}; their averages are over reduced words and for such distributions, one cannot use the theory of random matrix products. Cantrell and Pollicot prove CLTs and continuity of dynamical quantities when counting closed geodesics on variable negative curvature surfaces \cite{cantrell2022comparison}, but their results are different from ours as they also consider different distributions, averaging over words in $\pi_1(\Sigma)$ of length at most $n$. Nevertheless, one should expect at a high level at least that these various distributions on $\pi_1(\Sigma)$ should at least roughly produce similar behavior; by work of Kesten in \cite{kesten1959symmetric} a ``typical'' word of length $n$ corresponds to a homotopy class with algebraic length close to $\lambda n$ for $0<\lambda < 1$, for large $n$. (Note that $\lambda > 0$ requires non-amenability of the group in question; in particular this is true for our purposes but not true in total generality.)

\begin{remm}
Random walks have also been used spectacularly in a number of applications to geometric rigidity and classification theorems; for some examples, see the classification of stationary measures \cite{benoist2011mesures} and Poincar\'{e} recurrence for random walks on homogeneous spaces \cite{eskin2004recurrence}, and the classification of $\SL_2(\R)$-orbit closures on the moduli space of quadratic differentials \cite{eskin2015isolation}. A central feature in the proofs of these results are \textit{Margulis functions}, which informally characterize how much time a random walk on a space $X$ avoids a compact subset $E \subset X$. In future work, we plan to study how these uniform large deviation estimates for random matrix products are stably reflected in various avoidance principles in homogeneous dynamics.
\end{remm}

\subsection{Outline of the paper}

In \Cref{prelimsec}, we introduce the basic objects at play, such as the spaces of probability measures satisfying certain moment conditions, and the associated logarithmic/semi-logarithmic Wasserstein topologies. We demonstrate various fundamental properties of these spaces; many of these are standard but \Cref{wedgecontinuity} seems fundamentally new, and is a crucial technical result which allows us to establish a uniform separation between the first and second Lyapunov exponent, i.e. \Cref{unifgapbound}. The technical result is a certain continuity result for pushforwards under exterior powers; even the result we obtain is quantitatively weak and delicate, and it may be of independent interest.

In \Cref{cocyclethms}, we prove necessary large deviation estimates for martingales and Markov chains, which are uniform over compacta in the appropriate topologies. The proof of this result for martingales roughly follows from the carefully tracing dependencies in the pointwise version established by Benoist and Quint. The result for linear functionals on Markov chains could be proven using their approach and some technical innovations of ours, but for brevity, we elected to make use of new results of \cite{cai2024statisticalpropertiesmixingmarkov}. We introduce a quantitative formulation of ``almost invariant measures" (see \Cref{almostinvMeasures}) which may be new in order to establish these uniform results.

In \Cref{invsec}, we prove our main theorems regarding large deviation estimates for random matrix products, and technical lemmas concerning uniform estimates on the associated stationary measures. These technical lemmas regarding the stationary measure crucially rely on the aforementioned \Cref{wedgecontinuity}. These technical lemmas are necessary to extract large deviation estimates for matrix elements (i.e. $|\langle x, A_n^\mu y \rangle|$) but not for e.g. $\|A_n^\mu x\|$. We prove the regularity result \Cref{regularitythmintro} here; the proof synthesizes an approach used by Benoist and Quint and the technical lemmas proved earlier in the section. Finally, we prove our result concerning continuity of variance \Cref{varcont}; this requires establishing uniform bounds on a certain cocycle appearing in the variance formula \ref{varianceCocycleEq}.

In \Cref{applsec}, we apply our results to problems on random Schr\"odinger operators and random geodesics on hyperbolic surfaces. The main novelty in both settings comes from the abstract results which appear in previous sections.

\subsection{Acknowledgements} The first author thanks Lana Jitomirskaya for posing the question of localization for heavy tailed distributions, and for many enlightening discussions. The second author thanks his advisor Jesse Wolfson for his support and interest in this work. We thank Anton Gorodetski and Anthony Sanchez for their feedback on a previous version of this article. OH was supported in part by NSF Grants No. DMS-2052899 and DMS-2155211, and Simons Grant No. 896624. SR was supported in part by NSF Grants No. DMS-1944862 and DMS-2342135.

\section{Preliminaries and basic properties of the Wasserstein type topologies}\label{prelimsec}
	\subsection{Conventions}In this work, the letter P is unfortunately in quite high demand; we need notation for various spaces of probability measures, for the probability of an event, and for projectivization of vector spaces. We will use expressions like spaces $\mathcal{P}(X)$ to denote probability measures on $X$, $\Proj(\mathbb{K}^d)$ to denote the projectivization of $\mathbb{K}^d$, and $\P[A]$ to denote the probability of an event $A$. Note that e.g. $\P(\R^2)$ is not the real projective plane, but the real projective line.
	Similarly, we will study martingales, or more precisely martingale differences. We recall that given a filtration of $\sigma$-algebras $\mathcal{F}_n$ on a probability space, a sequence of variables $\vp_n$ is a \textit{martingale difference} if for all $n$, $\vp_n$ is $\mathcal{F}_n$ measurable and moreover $\E[\vp_n\,|\,\mathcal{F}_{n-1}] = 0$. For any random variable $X$, whether scalar or matrix valued, we let $\mathcal{L}(X)$ denotes its law.
	\subsection{Concave Wasserstein distances} Here we construct the metrics on (subspaces of) the space of probability measures, which we consider to be the appropriate topologies for working with heavy-tailed measures. We will make use of two elementary facts which are crucial:
	\begin{fact}
		Let $(M,d)$ be a metric space. If $f:[0,\infty) \rightarrow [0,\infty)$ is concave, strictly increasing, and satisfies $f(0) = 0$, then $d^f(x,y) := f(d(x,y))$ is also a metric on $M$.
	\end{fact}
	\begin{fact}\label{conccont}
		If $f:(a,\infty)\rightarrow \R$, for some $a \in \R$, is differentiable, $f'(x)$ is decreasing,  $\lim_{x\rightarrow+\infty} = +\infty$, and $\lim_{x\rightarrow +\infty} f'(x) = 0$, then there is a function $\tilde{f}:[0,\infty)\rightarrow [0,\infty)$ such that
		\begin{enumerate}
			\item $\tilde{f}(0) = 0$,
			\item $\tilde{f}(x) = f(x)$ for $x$ sufficiently large,
			\item $\tilde{f}$ is concave.
		\end{enumerate}
	\end{fact}
	Mostly to be able to give an explicit construction which we will take as definition, we prove the latter fact.
	\begin{proof}
		We will take $\tilde{f}$ to be of the form
		\begin{equation}\label{linint} \tilde{f}(x) =\begin{cases}
			\frac{f(x_0)}{x_0}x &\quad x \leq x_0\\
			f(x)&\quad x > x_0
		\end{cases}\end{equation}
		and it suffices to find an appropriate choice of $x_0$.
		Note that $f'(x)$ is actually continuous. Fix $b>0$ such that $\frac{f(b)}{b} > 0$. For all sufficiently large $x$, we obtain
		\[f'(x) \leq \frac{f(b)}{b}\]
		which readily implies
		\[f'(x)x \leq f'(x)(x-b) + f(b)\]
		The right hand side is an underestimate of $f(x)$ by the assumption that $f'(x)$ is decreasing. Hence, \begin{equation}\label{diffineq}
			f'(x) \leq \frac{f(x)}{x}
		\end{equation}
		If the infimum of all $x$ satisfying $(\ref{diffineq})$ is in $(a,\infty)$, we take $x_0$ to be this infimum, and in particular we will have $f'(x_0) = \frac{f(x_0)}{x_0}$. If this infimum is either negative, or equal to $a$, we just take $x_0 = \max\{a+1,1\}$. For such a choice of $x_0$, $\tilde{f}$ as defined in (\ref{linint}) is concave. (In the case where $f'(x_0) = \frac{f(x_0)}{x_0}$, $\tilde{f}$ is moreover differentiable.)
	\end{proof}
	Using these facts, we will be able to build various ``concave'' Wasserstein distances which nicely topologize spaces of distributions satisfying e.g. fractional or logarithmic moment conditions. The following families of functions satisfy the hypotheses of Fact \ref{conccont}:
	\begin{enumerate}
		\item $(\log x)^p$ for $p >0$,
		\item $e^{(\log x)^\delta}$ for $\delta \in (0,1)$,
        \item $x^p$ for $p \in (0,1)$.
	\end{enumerate}
	(Obviously, we do not need to use Fact \ref{conccont} for the third family.)
 
    Specifically, we define for $p \geq 1$ the family of functions:
	\[\log^{p\star}(t) := \begin{cases} \left(\frac{p}{e}\right)^pt&\text{ for }t < e^p\\ \log^p(t)&\text{ for }t \geq e^p\end{cases}\]
	These functions are concave, positive for positive inputs, and zero at zero. Hence we can define, for the original $d$, the logarithmic distances $d^{p\star} := d^{\log^{p\star}}$. For $\delta \in (0,1)$, by choosing an appropriate function corresponding to $e^{(\log x)^\delta}$, we can also construct semi-logarithmic distances. We let $\slog^{\delta}$ denote such a function, and carry out the corresponding construction, obtaining the semi-logarithmic distances $d^{\delta\star} := d^{\slog^\delta}$ for $\delta \in (0,1)$. Finally, this procedure allows one to define fractional distances $d^p$ for $p \in (0,1)$ by using the distance $d^p(g_1,g_2) := (d(g_1,g_2))^p$.
	
    By using these metrics, which give the space a very different large scale geometry, we are able to produce appropriate topologies for the study of distributions with heavy tails; the logarithmic distances will give rise to a topology appropriate for the study of measures with logarithmic moments, the semi-logarithmic distances give rise to a topology appropriate for the study of measures with semi-logarithmic moments, and the fractional power distances one appropriate for the study of measures with fractional moments. These topologies are precisely those generated by Wasserstein distances.
	
    Let $(G,d)$ be a locally compact metric group, and $\iota$ its identity. Most of this section generalizes to arbitrary metric spaces, but all our results are concerned with metric groups, and our main results concern the case $G=\mathrm{SL}_d(\mathbb{K})$ or $\GL_d(\mathbb{K})$ for $\mathbb{K}$ a local field. We recall briefly the definition of the Wasserstein (also known as Kantorovich--Rubinstein) distances.
	
	Given two distributions $\mu_1$ and $\mu_2$, both of which have a finite first moment in the sense that
	\[ \int d(g,\iota)d\mu_i(g) < \infty \]
	for $i=1,2$, the \textit{Wasserstein distance} $W^1$ between them is given by
	\[ W^1(\mu_1,\mu_2) := \inf_{\eta} \int d(g,h)d\eta(g,h) \]
	where the infimum is taken over all couplings $\eta$ of $\mu_1$ and $\mu_2$, i.e. distributions $\eta$ on $G^2$ such that $\eta$ has first marginal distribution $\mu_1$ and second marginal distribution $\mu_2$.
	
	Assuming higher moments, one can analogously define higher Wasserstein distances as follows for any $p \geq 1$:
    \begin{equation}
        W^p(\mu_1,\mu) := \inf_{\eta} \left[\int [d(g_1,g_2)]^p \,d\eta(g_1,g_2)\right]^{1/p}.
    \end{equation}
    However, for certain applications we also sought ``lower'' Wasserstein distances. In this case, rather than building something analogous at the level of the metric on distributions though, one just changes the metric.
	
	Hence, we define the \textit{logarithmic Wasserstein distances }$W^p_{\log}$ as follows:
	\[W^{p}_{\log}(\mu_1,\mu_2) := \inf_\eta \int d^{p\star}(g,h) d\eta(g,h)\]
    where the infimum is taken over all couplings $\eta$ of the measures $\mu_1$ and $\mu_2$.
	Obviously this is not well-defined for all distributions; we require a logarithmic moment of order $p$, i.e. \[\int_G d^{p\star}(g,\iota)d\mu_i(g) = \int_G \log^{p\star}(d(g,\iota))d\mu_i(g) < \infty. \]
	Similarly, we will define \textit{semi-logarithmic Wasserstein distances} for $0 < \delta < 1$ by using $d^{\delta \star}$; yielding the semi-logarithmic Wasserstein distances which we denote by $W^{\delta}_{\slog}$. The \textit{fractional Wasserstein distances} $W^p$, for $p \in (0,1)$ are also defined analogously. We note that the fractional Wasserstein distances have already been used fruitfully in the study of random walks on linear groups; Tall and Viana used these distances to study the modulus of continuity for Lyapunov exponents of compactly supported distributions on $\mathrm{GL}_2(\R)$ \cite{tall2020moduli}.
	
	We emphasize that with the logarithmic, semi-logarithmic, and fractional Wasserstein distances, we are just using $W^1$ with a different choice of metric on $G$. Having said this, we will now use $W^1$ going forward to denote the Wasserstein metric (on spaces of probability measures) associated to the original metric on the underlying space; this should not lead to too much ambiguity, as we will principally be concerned with the very explicit case where the space is $\mathrm{SL}_d(\mathbb{K})$ and the ``original'' distance is that induced by the operator norm.
	
	The following fact is well-known, and central to our work.
	\begin{prop}\label{momentmetric}

    Let $(G,d)$ be a metric group, and $W^p$ an associated Wasserstein distance for $p \geq 1$. Then 
		$W^p(\mu_n,\mu) \rightarrow 0$ if and only if $\mu_n \rightarrow \mu$ weakly and moreover the $p$-th moments of $\mu_n$ converge to that of $\mu$, i.e.
		\[ \int [d(g,\iota)]^pd\mu_n(g) \rightarrow \int [d(g,\iota)]^pd\mu(g)\]
\end{prop}

This has as an immediate corollary that convergence in a logarithmic (resp. semi-logarithmic, fractional) Wasserstein topology is equivalent to weak convergence and convergence of the corresponding logarithmic (resp. semi-logarithmic, fractional) moment.

Another important and well-known fact regarding the Wasserstein metric is \emph{Kantorovich duality}, which gives us the formulae:
\begin{equation}
		W^{p}_{\log}(\mu,\mu') = \sup\left\{\int f\,d(\mu-\mu')\,:\, \sup_{g\neq g'\in G}\frac{f(g)-f(g')}{\log^{p\star}(d(g,g'))} \leq 1 \right\}
	\end{equation}
	\begin{equation}
		W^{\delta}_{\slog}(\mu,\mu') = \sup\left\{\int f\,d(\mu-\mu')\,:\, \sup_{g\neq g'\in G}\frac{f(g)-f(g')}{\slog^{\delta}(d(g,g'))} \leq 1 \right\}
\end{equation}
\begin{equation}
    W^1(\mu,\mu') = \sup\left\{ \int f d(\mu-\mu')\,:\, \sup_{g \neq g' \in G} \frac{f(g) - f(g')}{d(g,g')} \leq 1 \right\}
\end{equation}
i.e. $W^1$ with respect to any distance has a variational characterization with respect to the functions into $\R$ which are 1-Lipschitz with respect to that distance.
	\subsection{Random matrix products}\label{rmp}
	Throughout the rest of the paper, we will mostly study the specific case where $G=\mathrm{SL}_d(\mathbb{K})$ equipped with the distance
	\[ d(A,B) := \|A-B\|\]
    where $\|A\|$ is the spectral norm. (We fix a norm to be definite, but the choice of norm does not matter; over local fields, all norms on a finite dimensional vector space are equivalent.)
\begin{remm}\label{glvssl}
Our results for $G=\SL_d(\KK)$ generalize to $G=\mathrm{GL}_d(\mathbb{K})$, with two very important caveats: if $G = \GL_d(\KK)$ one must then work with the distance $d(A,B) = \max\{\|A-B\|,\|A^{-1}-B^{-1}\|\}$ instead in order to control the least singular value. Moreover, our results which we prove for $\SL_2(\KK)$ specifically cannot be straightforwardly generalized to $\GL_2(\KK)$ by our methods. To simplify exposition, we treat the $\SL_d(\KK)$ case explicitly.
    
\end{remm}

    In particular, we are for the most part interested in studying stability or continuity of various statistical quantities associated to certain types of random walks on $\mathrm{SL}_d(\mathbb{K})$, specifically those associated to random matrix products. We recall the following basic facts regarding random matrix products.
    
    Given any distribution $\mu$ on $\mathrm{SL}_d(\mathbb{K})$, a product of $n$ matrices independently distributed with law $\mu$ has distribution $\mu^{\ast n}$ with $\mu^{\ast 1} = \mu$ and
    \begin{equation}
        \int f(g) d\mu^{\ast n}(g) = \int f(g_1g_2) \,d\mu(g_1)\,d\mu^{\ast (n-1)}(g_2)
    \end{equation}
    for $n > 1$.

    To formulate things more succinctly, here and throughout we let $A_n^\mu$ be a $\mu^{\ast n}$ distributed $\mathrm{SL}_d(\mathbb{K})$-valued random variable. Under a very lenient moment assumption, specifically
    \begin{equation}\label{furstmoment}
		\int \max\{\log\|A\|,0\}d\mu(A) < \infty,
	\end{equation}
one can define Lyapunov exponents as follows:
	\begin{equation} \lambda_1(\mu):=\lim_{n\rightarrow \infty} \frac{1}{n}\E[\log\|A_n^\mu\|].\end{equation}
	For $1 < m \leq d$, we define them inductively by
	\begin{equation}
		\sum_{k=1}^m \lambda_k(\mu) = \lim_{n\rightarrow \infty} \frac{1}{n} \E[\log \| \wedge^m A_n^\mu\|]
	\end{equation}
    with $\wedge^m$ denoting the $m$-th exterior power. (Later in this paper, these quantities are always finite, but in the generality specified here, these quantities can be $-\infty$.) Note that existence of these quantities is a straightforward consequence of Fekete's lemma, and that these values give the asymptotic behavior of the singular values of random matrix products, in the sense that ``typically'' $s_m(A_n^\mu)$ is close to $e^{\lambda_m n}$, where $s_1(A_n^\mu) \leq \cdots \leq s_d(A_n^\mu)$ are the singular values of $A_n^\mu$.

    These quantities encode the asymptotic behavior of random matrix products of $\mu$ in a form made precise by Osceledets' theorem. For our purposes, we are concerned almost entirely with $\lambda_1$ and $\lambda_2$; $\lambda_1$ encodes the almost sure asymptotic behavior of random matrix products of the form $A_N = M_N\cdots M_1$ with $M_n$ all i.i.d. with law $\mu$, and control of $\lambda_1-\lambda_2$ is important in getting uniform estimates on asymptotic almost sure behavior. One has the following almost surely:
	\begin{equation}
		\frac{1}{N}\log\|A_N^\mu\| \rightarrow \lambda_1(\mu),
	\end{equation}
	and for any fixed non-zero $x \in \R^d$:
	\begin{equation}
		\frac{1}{N} \log \|A_N^\mu x\| \rightarrow \lambda_1(\mu).
	\end{equation}
	These facts are now classical, first proven by Furstenberg; they also follow easily from Kingman's subadditive theorem (proved after the original proof of Furstenberg). In particular, in order to prove it via Kingman's subadditive theorem one models the random matrix product as a matrix valued functional of a Bernoulli shift. It is not hard to see that (\ref{furstmoment}) is implied by
	\[ \int d^{1\star}(g,\iota)d\mu(g) < \infty\]
	These notions are not unique to the study of random walks of this type; the behavior of $\log\|A_N x\|$ is best studied in terms of the cocycle $\Phi(g,x): \SL_d(\mathbb{K})\times \Proj (\mathbb{K}^d) \to \KK$ defined by $\Phi(g,x) = \log\frac{\|gx\|}{\|x\|}$; this is called the \textit{log-norm cocycle} and has been extensively studied with $g$ given by various underlying dynamics. These basic results hold for any underlying ergodic dynamics so long as the moment assumption holds. However, the finer features of the statistical behavior of these random matrix products are very sensitive to the underlying dynamics. 
	
	In this work we focus exclusively on the case of random i.i.d. matrix products, which essentially comes from the underlying dynamics being a Bernoulli shift (with possibly infinite alphabet). Our results all concern the continuity and/or stability of asymptotic bounds or behavior for such products in the topologies we introduce. 
	These Lyapunov exponents are particularly well-behaved under certain dynamical assumptions on the random walk.
 
	\begin{deffo}\label{nddcond} Let $\Gamma_\mu$ denote the subgroup of $\SL_d(\mathbb{K})$ generated by the support of a measure $\mu$. We say $\Gamma_\mu$ is \emph{strongly irreducible} if there is no finite collection of linear subspaces $V_i \subset \mathbb{K}^d$ such that $g[\cup V_i] \subset \cup V_i$ for all $g \in \Gamma_\mu$, i.e. there is no finite union of subspaces invariant under the action of all $g \in \Gamma_\mu$. We say $\Gamma_\mu$ is \emph{contracting} if there is a sequence $g_n$ of elements in $\Gamma_\mu$ such that $\frac{g_n}{\|g_n\|}$ converge to a rank one matrix.
	\end{deffo}
	A result of Furstenberg and Kifer \cite{FurstenbergKifer} in particular implies that if $\Gamma_\mu$ is strongly irreducible, contracting, and (\ref{furstmoment}) holds, then $W^p_{\log}(\mu_n,\mu) \rightarrow 0$ (henceforth written $\mu_n \overset{W^p_{\log}}{\rightarrow}\mu$) implies that $\lambda_1(\mu_n) \rightarrow \lambda_1(\mu)$. (The $\mu_n$ do not need to satisfy the dynamical assumptions.)
	
	While we cannot hope for continuity of the top Lyapunov exponent $\lambda_1$ outside the locus of measures satisfying certain constraints such as strong irreducibility and contraction, at least not in such a coarse topology (c.f. \cite{avila2023continuitylyapunovexponentsrandom,bocker2017continuity}), upper semi-continuity was shown for the $W^1$ Wasserstein topology on $\SL_2(\R)$ by S\'anchez and Viana:
	
	\begin{thm}[\cite{sánchez2020lyapunovexponentsprobabilitydistributions}]
		On the space of distributions $\mu$ on satisfying
		\[ \int d(g,\iota) \,d\mu(g) < \infty\]
		equipped with the $W^1$ metric, $\mu \mapsto \lambda_1(\mu)$ defines an upper semi-continuous function.
		
	\end{thm}
	
	Their proof readily generalizes, allowing us to prove the following; we emphasize that our proof for our various concave Wasserstein topologies is essentially that of S\'anchez and Viana for $W^1$.
	\begin{thm}\label{usclyap}
		Fix $p>0$. The function $\mu\mapsto \lambda_1(\mu)$ is upper semi-continuous in the $W^1_{\log}(\SL_d(\KK))$ and hence in all the finer topologies.
	\end{thm}
	
	We will reproduce the proof below; it introduces many useful notions in a concrete setting, and the argument of S\'anchez and Viana is also necessary for technical reasons in our results on general cocycles. 
    
    The group $\SL_d(\KK)$ naturally acts linearly on the associated projective space $\P(\KK^d)$, and so we can develop a notion of measure invariance via convolution of measures on $\SL_d(\KK)$ and $\P(\KK^d)$. Given a measure $\mu$ on $\SL_d(\KK)$, we say a measure $\nu$ on $\P(\KK^d)$ is \emph{$\mu$-invariant} (or \emph{$\mu$-stationary}) if $\mu * \nu = \nu$, where this convolution is a measure on the projective space defined by
    \[\int_{\P(\KK^d)} f(x) d(\mu \ast \nu) = \int_{\SL_d(\KK)\times\P(\KK^d)} f(Mx) \,d\mu(M)\,d\nu(x) .\]
    Such measures exist under our dynamical assumptions; see the beginning of Section \ref{cocyclethms} for more details in greater generality.
	\begin{fact}
		Given any $\mu$ satisfying (\ref{furstmoment}),
		\begin{equation}
			\lambda_1(\mu):= \max_\nu \int \frac{\|A x\|}{\|x\|}d\mu(A)d\nu(x)
		\end{equation}
		where the maximum is taken over all $\mu$-invariant measures $\nu$ on projective space.
	\end{fact}
	\begin{fact}
		If $\mu_k \rightarrow \mu$ weakly, and $\nu_k$ are $\mu_k$-invariant measures on the associated projective space, any weak limit point of the sequence $\nu_k$ is $\mu$ invariant.
	\end{fact}
	We use these facts to demonstrate the claimed upper semi-continuity by the strategy pursued in \cite{sánchez2020lyapunovexponentsprobabilitydistributions}. In order to do this, we first establish the following:

\begin{prop}\label{pthLogDominatingNormCocycle}
    Fix some $p >1$. Away from a compact set $K' \subset \SL_d(\KK)$, the $p$-th logarithmic Wasserstein distance dominates the log-norm cocycle $\Phi$, i.e.
    $$|\Phi(g,[v])| \leq C d^{p\star}(g,\iota)$$
    for some constant $C$.
\end{prop}
\begin{proof}
    Let us first give a general bound on the log-norm cocycle. Since the function $\log : [1, \infty) \to \R$ is $1$-Lipschitz and $\|g\| \geq 1$ for all $g \in \SL_d(\KK)$, we have
    \begin{equation}\label{NormDistanceBound}
        |\Phi(g,[v])| = \left|\log \frac{\|g v\|}{\|v\|}\right| \leq \log\|g\| \leq \big|\|g\| - \|\iota\|\big| \leq \|g - \iota\|=d(g,\iota).    \end{equation}
    To prove the claim, we split into two main cases. We first deal with the case of small distances, that is, we assume $d(g,\iota) < e^p$. In this case, $d^{p\star}(g,\iota) = (p/e)^p d(g,\iota)$. When $p \geq e$, we have
    $$d(g,\iota) \leq \left(\frac{p}{e}\right)^p d(g,\iota) = d^{p\star}(g,\iota),$$
    When $1\leq p < e$, we trivially have
    $$d(g,\iota) \leq \left(\frac{e}{p}\right)^p d^{p\star}(g,\iota),$$
    proving the desired domination in this case.

    We now deal with the case of large distances, that is, we assume $d(g,\iota) \geq e^p$. In this case, $d^{p\star}(g,\iota) = \log^p(d(g,\iota))$. We do not deal with the compact subset of bounded norm matrices $g \in\SL_d(\KK)$ where
    \begin{equation}\label{enlargedpthLogCompact}
        K' = \{e^p \leq d(g,\iota) \leq \|g\|\},
    \end{equation}
    for in this range the domination does not hold. We then further assume that $d(g,\iota) > \|g\|$ for sufficiently large $\|g\|$. If this is this case, then taking $\log$ of both sides yields
    $$\log\|g\| < \log(d(g,\iota)),$$
    and because $x < x^p$ for $p>1$ and $x \geq 1$, we have $\log(d(g,\iota)) < \log^p(d(g,\iota)) = d^{p\star}(g,\iota)$. By \ref{NormDistanceBound}, we have thus shown
    $$|\Phi(g,[v])| \leq d^{p\star} (g,\iota),$$
    proving the claim.
\end{proof}

We can now prove upper semi-continuity of the Lyapunov map for the $p$-th logarithmic Wasserstein probability space:

	\begin{proof}[Proof of \Cref{usclyap}]
                By the definition of upper semicontinuity, need to show that given a convergent sequence of probability measures $\mu_k \xrightarrow[]{W^p_{\log}} \mu$, the corresponding sequence of Lyapunov exponents satisfies $\lambda_1(\mu_k)\leq \lambda_1(\mu)$. We will need to use facts about convergence in the $p$-th logarithmic Wasserstein topology. 

    For each probability measure $\mu_k$ on $\SL_d(\KK)$, let $\nu_k$ denote a $\mu_k$-stationary probability measure on $\Proj(\KK^n)$ achieving the maximum given by \cite[Proposition 6.7]{viana2014lectures}:
    $$\lambda_1(\mu_k) = \int \Phi d\mu_k d\nu_k.$$ By compactness of $\Proj(\KK^n)$, we can pass to a subsequence $\nu_{k_j}$ which converges in the weak-$*$ topology to a $\mu$-stationary measure $\nu$. For any $\epsilon > 0$, we want to show that there exists some $k_0 \in \N$ so that for all $k > k_0$,
    $$\left|\int \Phi d\mu_k d\nu_k - \int \Phi d\mu d\nu\right| < \varepsilon.$$
    Since the $p$-th logarithmic moment of $\mu$ is finite, there exists some compact subset $K_1\subset \SL_d(\KK)$ so that
    \begin{equation}
        \int_{K_1} d^{p\star}(g,\iota)d\mu(g) < \frac{\varepsilon}{36}.
    \end{equation}
    Moreover, since $W^p_{\log}(\mu_k , \mu) \to 0$, we have that there is some $R' >0$ and $k' \in \N$ so that for all $k > k'$, we have
    \begin{equation}
        \int_{d^{p\star}(g,\iota) > R'} d^{p\star}(g,\iota)d\mu_k(g) < \frac{\varepsilon}{36}.
    \end{equation}
    Let $R > 0$ be big enough so that $B_\iota (R')\cup K_1 \subset B_\iota (R)$ and define the compact set $K = \overline{B_\iota(R)} \cup K'$, where $K'$ is given in \ref{enlargedpthLogCompact}.
    We observe that 
        \begin{align*}
            &\left|\int \Phi d\mu_k d\nu_k - \int \Phi d\mu d\nu\right| \\
        &\leq  \left|\int_{K\times \Proj(\KK^n)} \Phi d\mu_k d\nu_k - \int_{K\times \Proj(\KK^n)} \Phi d\mu d\nu\right|  +\left|\int_{K^c\times \Proj(\KK^n)} \Phi d\mu_k d\nu_k \right| + \left| \int_{K^c\times \Proj(\KK^n)} \Phi d\mu d\nu\right|. 
        \end{align*}
    On the last two terms, the inequality $|\Phi(g,[v])| \leq d^{p\star}(g,\iota)$ hold by \Cref{pthLogDominatingNormCocycle}, and so they are each bounded by $\varepsilon/3$. The first term can be readily bounded by $\varepsilon/3$ in the same manner that S\'{a}nchez--Viana carries out \textit{mutatis mutandis} (using compactness and the Ursohyn function technique), therefore we can conclude
    $$\left|\int \Phi d\mu_k d\nu_k - \int \Phi d\mu d\nu\right| < \varepsilon.$$
    From this, we conclude that
    $$\lambda_1(\mu_k) = \int \Phi d\mu_kd\nu_k \longrightarrow \int \Phi d\mu d\nu \leq \lambda_1(\mu),$$
    thereby proving upper semicontinuity of the first Lyapunov exponent in the $p$-th logarithmic Wasserstein topology.
	\end{proof}
Continuity of another map is also important. Recall that $\wedge^2 A$ denotes the second exterior power of a matrix $A$. Fixing a choice of basis for $\wedge^2 \KK^n$, $\wedge^2$ induces a map from $\SL_d(\mathbb{K})$ to $\SL_{\binom{d}{2}}(\KK)$. Using the bound $\|\wedge^2A\| \leq \|A\|^2$, it is obvious that
\begin{enumerate}
    \item If $\mu \in \mathcal{P}^{p}_{\log}(\SL_d(\KK)$, then $\wedge^2_\ast \mu \in \mathcal{P}^{p}_{\log}(\SL_{\binom{d}{2}}(\KK))$.
    \item If $\mu \in \mathcal{P}^{\delta}_{\slog}$, then $\wedge^2_\ast \mu \in \mathcal{P}^{\delta'}_{\slog}(\SL_{\binom{d}{2}}(\KK))$ for any $
\delta' < \delta$.
    \item If $\mu \in \mathcal{P}^p(\SL_d(\KK))$, then $\wedge^2_\ast \mu \in \mathcal{P}^{p/2}(\SL_{\binom{d}{2}}(\KK))$.
\end{enumerate}
However, continuity of the mapping $\mu \mapsto \wedge^2_\ast \mu$ is considerably more delicate. In order to establish lower semi-continuity of the quantity $\lambda_1 - \lambda_2$ later, it is useful to at the very least establish continuity from various ``higher'' Wasserstein spaces into $\mathcal{P}^1_{\log}$, as $\lambda_2$ is upper semi-continuous with respect to this topology.
\begin{thm}\label{wedgecontinuity}
    For any $\delta \in (0,1)$, the map $\mu \mapsto \wedge^2_\ast \mu$ from $\mathcal{P}^\delta_{\slog}(\SL_d(\KK))$ to $\mathcal{P}^1_{\log}(\SL_{\binom{d}{2}}(\KK))$ is continuous.
\end{thm}
\begin{proof}
    Take $\mu'$ such that $W^\delta_{\slog}(\mu',\mu) < \ve$ for some small $\ve$; in particular we can take $A,B$ to be $\mu$ and $\mu'$ distributed variables respectively such that $\E[\slog^\delta(\|A-B\|)] < 2\ve$. In order to bound $\E[\log^{1\star}(\|\wedge^2A - \wedge^2 B\|)]$, we bound the variables $\E[\log^{1\star}(\|\wedge^2A - \wedge^2 B\|)1_{\{\|A-B\| > \ve^{-1/3}\}}]$ and $\E[\log^{1\star}(\|\wedge^2A - \wedge^2 B\|)1_{\{\|A-B\| \leq  \ve^{-1/3}\}}]$. We will use the straightforward bound
    \begin{equation}
        \|\wedge^2A - \wedge^2B\| \leq \|A\|\cdot \|A-B\| + \|B\|\cdot\|A-B\|
    \end{equation}
    and the concavity of $\log^{1\star}$ and $\slog^\delta$ throughout. Recall also that for $\ve$ sufficiently small (depending on $A$), $\E[\slog^\delta(\|B\|)] \leq \E[\slog^\delta(\|A\|)] + 1$. Finally, we will use the fact that $\log^{1\star}(xy) \leq C(\log^{1\star}(x) + \log^{1\star}(y))$ freely below.

    Combining all these, we get
    \begin{align*}
        \E[\log^{1\star}(\|\wedge^2A - \wedge^2 B\|)1_{\{\|A-B\| > \ve^{1/3}\}}] &\leq (\E[\log^{1\star}(\|\wedge^2A - \wedge^2 B\|)^2]\P[\|A-B\| > \ve^{1/3}])^{1/2}\\
        &\leq(\E[\log^{1\star}((\|A\|+\|B\|)\ve^{-1/3})^2]\P[\|A-B\| > \ve^{1/3}])^{1/2}\\
        &\leq C(\E[(\log^{1\star}(\|A\|+\|B\|)-\frac{1}{3}\log(\ve))^2]\P[\|A-B\| > \ve^{1/3}])^{1/2}\\
        &\leq  C_\mu \ve^{1/2} \log(1/\ve) (\slog^\delta(\ve^{1/3}))^{-1/2} \underset{\ve \rightarrow 0}{\rightarrow} 0
    \end{align*}
    using Cauchy--Schwartz and Chebyshev in the last step.

    To bound the other piece, we use the fact that
    \begin{equation}\label{thinggg}
    \E[\log^{1\star}(\|\wedge^2A - \wedge^2 B\|)1_{\{\|A-B\| \leq \ve^{1/3}\}}] \leq \sum_{k=1}^\infty \P[\log^{1\star}((\|A\|+\|B\|)\ve^{1/3}) \geq k\ve^{1/3}] + \ve^{1/3}
    \end{equation}
    By Chebyshev, the right hand side is bounded by:
    \[
    \ve^{1/3} + C_\mu \ve^{1/3}\sum_{k=1}^\infty e^{-(k\ve^{1/3})^\delta}
    \]
    so that in particular it suffices to show that the infinite sum is $o(\ve^{-1/3})$ as $\ve \rightarrow 0$. However, this is an immediate consequence of the fact that the largest $k$ such that $e^{-(k\ve^{1/3})^\delta} > k^{-2}$ is $o(\ve^{-1/3})$, so that in particular the right hand side of (\ref{thinggg}) goes to zero as $\ve \rightarrow 0$.
\end{proof}

This theorem, combined with Proposition \ref{usclyap}, yields \Cref{unifgapbound} and \Cref{uscsum} as corollaries.

\begin{remm}\label{wedgecontrmk}
    Continuity of $\mu \mapsto \Lambda^2_\ast \mu$ suffices for our purposes; note that at the cost of somewhat messier estimates, one can use this approach to obtain continuity of all exterior power pushforwards as maps from the semi-logarithmic Wasserstein spaces $W^\delta_{\slog}(\SL_d(\KK))$ to $W^{1}_{\log}(\SL_{\binom{d}{k}}(\KK))$, or any logarithmic Wasserstein space. In particular, \Cref{uscsum} can be generalized to hold for $\lambda_1(\mu) + \cdots + \lambda_k(\mu)$ for any $k \leq d$.
    
    Moreover, we mention that for $\SL_d(\KK)$ cocycles, it is trivial that pushforward by $\wedge^d$ is continuous regardless of moment assumptions; $\wedge^d$ is the determinant. In particular $\lambda_1(\mu)+\lambda_2(\mu)$ is an upper semicontinuous function of $\mu \in W^p_{\log}(\SL_2(\KK))$ for any $p \geq 1$; we emphasize this trivial fact because it explains the difference in results obtained for $\SL_2(\KK)$-cocycles and those obtained for $d>2$.
\end{remm}
\section{Uniform large deviation estimates for martingales and Markov chains}\label{cocyclethms}
In this section, we will obtain most of the estimates needed to prove our large deviation results, and one of our results will hold fairly generally. If $G$ is a locally compact topological group, and $X$ is a compact metrizable $G$-space, (by $G$-space we just mean equipped with a continuous $G$-action), then any function $\Phi: G\times X \rightarrow \KK$ is called a \textit{cocycle} if it satisfies the identity
\begin{equation}
	\Phi(g_1g_2,x) = \Phi(g_1,g_2x) + \Phi(g_2,x).
\end{equation}
In particular, one can study random matrix products via the log-norm cocycle $\Phi(g,x)= \log \frac{\|gx\|}{\|x\|}$ for $G = \GL_d(\mathbb{K})$ and $X = \P(\mathbb{K}^d)$. (We abuse notation, letting $x$ represent both a point in $\P(\mathbb{K}^d))$ and also some vector in $\mathbb{K}^d$ pointing in said direction.) Indeed, if $g_k$ are random matrices, i.i.d. with common law $\mu$, then the random product has distribution given by a convolution power $\mu^{\ast n}$ as defined above. This framework also allows the study of random products in more general groups, and some of our results will hold in this generality. 

One can also convolve measures on the group $G$ with measures on $X$; for $\mu$ a measure on $G$ and $\nu$ a measure on $X$, we define the measure $\mu \ast \nu$ on $X$ by

\begin{equation}
	\int_X f(x)\,d(\mu \ast \nu)(x) = \int_{G\times X} f(gx)\,d\mu(g)\,d\nu(x)\quad \text{ for all } f \in C(X)
\end{equation}

This notion is largely important because of the importance of invariant measures:
\begin{deffo}
	Given a measure $\mu$ on $G$, we call a measure $\nu$ on $X$ \textit{$\mu$-stationary}, or \textit{$\mu$-invariant}, if $\mu \ast \nu = \nu$. For any $\mu$, we let $\mathrm{Inv}(\mu)$ denote the set of $\mu$-invariant measures.
\end{deffo}
It is classical that any $\mu$ admits at least one stationary measure. For many of our specific results regarding random matrix products, we will assume hypotheses which guarantee there is exactly one --- this is one of the reasons for our assumption of strong irreducibility and contraction in our main results on random walks on linear groups.

Given these stationary measures, we can define a notion necessary to formulate our main result regarding cocycles:
\begin{deffo}\label{def:uniqueAvg}
    We say $\Phi \in C(X)$ has a \textit{unique average} (with respect to $\mu$) if the set of values
    \[\left\{\int \Phi \,d\nu\,:\, \nu \in \mathrm{Inv}(\mu)\right\}\]
    has a single element.
\end{deffo}

Clearly if $\mathrm{Inv}(\mu)$ has only one element, then any $\Phi \in C(X)$ has a unique average. In particular, we introduce the notion of unique average to prove results in the highest generality, but in our applications the stronger assumption that there is a unique invariant measure will always hold (strong irreducibility of $\mu$ guarantees this). Our main results in this section are large deviation theorems for cocycles which have a unique average. For the first family of results, we will treat cocycles associated to random matrix products with a strongly irreducible and contracting measure, whereas the last result is more general.

The central idea towards proving these results is to replace the cocycle $\Phi$ with
\begin{equation}\label{mgmarkovdecomp}\Phi'(g,x) := \Phi(g,x) - \int_{g\in G}\Phi(g,x)\,d\mu(g).\end{equation}

If one lets $g_k$ be i.i.d. distributed for $k \in \N$ with distribution $\mu$ and defines $h_n := g_n\cdots g_2g_1$, then $h_n$ has distribution $\mu^{\ast n}$ and $\Phi'(h_n,x)$ is a martingale with respect to the filtration $\mathcal{F}_n := \sigma(h_1,\dots,h_n)$. To get uniform large deviation estimates for $\Phi'(g,x)$ in $\mu^{\ast n}$, it suffices to get uniform large deviation estimates for martingales. We do this for the case of logarithmic moments in Theorem \ref{explicitmgpolylde}; in the semi-logarithmic and fractional/polynomial moment such results are well-known via the standard Chernoff bound method used to obtain e.g. Azuma's bound for bounded martingales.

The function $\int_{g\in G}\Phi(g,x)\,d\mu(g)$ has no dependence on $g$, and is continuous in $x \in X$. We will derive exponential large deviation bounds for this quantity in all regimes in Lemma \ref{unifblln}; in particular (by the cocycle property of $\Phi$), it is in fact a Birkhoff average corresponding to a continuous functional with respect to a certain Markov chain.

\subsection{Uniform large deviations for martingale differences}
We fix a probability space $(\Omega, \mathcal{F},\P)$ and a filtration of $\sigma$-algebras $\mathcal{F}_0 \subset \mathcal{F}_1 \subset \cdots$ on it. We recall the classical notion of a martingale difference:
\begin{deffo}
	A sequence of random variables $\vp_1,\vp_2,\dots$ (finite or infinite) is a \textit{martingale difference} if $\vp_n$ is $\mathcal{F}_n$-measurable and $\E[\vp_n\,|\,\mathcal{F}_{n-1}] = 0$.
\end{deffo}

Given a martingale difference, either infinite or of length at least $n$, one can study the variable 
\begin{equation}\label{summ}
	S_n := \sum_{k=1}^n \vp_n
\end{equation}
The sequence $S_1,\dots, S_n$ in fact is a martingale, and importantly $\E[S_k] = 0$. A central concern in this paper is controlling the probability of deviations from 0 for these variables. A crucial result is \cite[Theorem 2.2]{bqclt}, formulated below:

\begin{thm}[\cite{bqclt}]\label{mgpolylde}
Let $p > 1$ and $(\vp_n)_{n\geq 1}$ be a martingale difference, with $S_n$ defined as in (\ref{summ}). Assume there is $\vp \in L^p(\Omega)$ so that for all $n \geq 1$ and $t>0$ we have
\begin{equation}\label{phidom}\P[|\vp_n| > t\,|\,\mathcal{F}_{n-1}] \leq \P[\vp > t]
\end{equation} almost surely.
Then there exist $C_n = C_n(p,\ve,\vp)$ such that for all $n \geq 1$ and $\ve > 0$ satisfying
\begin{enumerate}
	\item $\P[|S_n| > n\ve] < C_n$
	\item $\sum_{n\geq 1} n^{p-2} C_n < \infty$
\end{enumerate}
\end{thm}
This result is not uniform; the constants depend on $\vp$, and in particular on the moments and the tail of $\vp$. (The fact that moment control alone does not suffice is demonstrated by an example in \cite[Remark 2.3]{bqclt}.)

To prove the result above, Benoist and Quint essentially proved the following; the proof we present is essentially theirs.
\begin{thm}\label{explicitmgpolylde}
Fix $p>1$ and $\ve > 0$ and define $\gamma := \frac{3p+1}{4p}$. Let $\vp_n$ be a martingale difference satisfying the uniform moment and tail conditions
	\begin{equation}\label{unimomentcontrol}
 N_1:= \sup_n \E[\vp_n^p\,|\,\mathcal{F}_{n-1}]< \infty\end{equation}
	and
	\begin{equation}\label{unitailcontrol}N_2:= \sup_{n\in\N} \sup\{ t \in \R\,: \E[\E[\vp_n^\gamma1_{\{\vp_n^\gamma > t\}}\,|\,\mathcal{F}_{n-1}]] > \ve/3 \} <\infty. \end{equation}
    Then there are constants $C_n:=C_n(p,\ve,N_1,N_2)$ such that $\sum_{n\geq 1}C_nn^{p-2}$ and
    \[\P[|S_n| > n\ve] \leq C_n]\]
\end{thm}
\begin{proof}
Recall that we are interested in proving large deviation estimates for a mean zero martingale $S_n = \sum_{1 \leq k \leq n} \vp_k$, under the following assumptions:
\begin{equation}\label{unimomentcontrolapp}N_1:= \sup_n \E[\vp_n^p\,|\,\mathcal{F}_{n-1}]< \infty\end{equation}
	and
	\begin{equation}\label{unitailcontrolapp}N_2:= \sup_{n\in\N} \sup\{ t \in \R\,: \E[\E[\vp_n^\gamma1_{\{\vp_n^\gamma > t\}}\,|\,\mathcal{F}_{n-1}]] > \ve/3 \} <\infty. \end{equation}
    Then there are constants $C_n:=C_n(p,\ve,N_1,N_2)$ such that $\sum_{n\geq 1}C_nn^{p-2}$ and
    \[\P[|S_n| > n\ve] \leq C_n].\]
We introduce the truncated variables
\[\vp_{n,k} = \vp_k1_{\{\vp_k< n^\gamma\}},\]
partial sums of these variables
\[T_n =\sum_{1\leq k \leq n} \vp_{n,k},\]
and also the following:
\[\overline{\vp_{n,k}} = \vp_{n,k} - \E[\vp_{n,k}\,|\,\mathcal{F}_{k-1}],\quad \overline{T}_n = \sum_{1\leq k \leq n} \overline{\vp_{n,k}}.\]
One obtains easily
\begin{equation}
	S_n = \overline{T}_n +(T_n - \overline{T}_n) + \sum_{1\leq k \leq n} \vp_k1_{\{|\vp_k| > n^\gamma\}}.
\end{equation}
The event $|S_n|> n\ve$ is thus necessarily (for large $n$) a subevent of the union of the following four events:
\begin{align*}
	A_{1,n}&:=  \{\text{there exists }k\leq n\text{ such that }|\vp_k| > \frac{\ve n}{3}\}\\
	A_{2,n}&:= \{\text{there exist }k_1 < k_2 \leq n\text{ such that }|\vp_{k_1}| > n^\gamma,\, |\vp_{k_2}| > n^\gamma\}\\
	A_{3,n}&:= \{ |T_n -\overline{T}_n| > \frac{\ve n}{3}\}\\
	A_{4,n}&:= \{ |\overline{T}_n| > \frac{\ve n }{3}\}
\end{align*}
It suffices to find $C_{i,n}$ so that $\P[A_{i,n}] \leq C_{i,n}$ and which satisfy the summability condition \begin{equation}\label{summc}
	\sum_{n\geq 1} n^{p-2}C_{i,n} < \infty
\end{equation}
and can be taken uniformly over all martingale differences satisfying uniform moment and tail bounds. In particular, we will show that $\P[A_{1,n}]$, $\P[A_{2,n}]$ and $\P[A_{4,n}]$ can be bounded in terms of moments of the martingale increments $\vp_k$ for $1 \leq k \leq n$; $\P[A_{3,n}]$ can be bounded in terms of quantities related to the tails of the increments with $1 \leq k \leq n$.

We start with the bounds on $\P[A_{1,n}]$. We have
\begin{align*}
    \P[A_{1,n}] &\leq \sum_{1 \leq k \leq n} \P[\vp_k > \ve n/3]\\
    &=\sum_{1 \leq k \leq n} \E[\P[\vp_k > \ve n/3\,|\,\mathcal{F}_{n-1}]]\\
    &\leq (\ve n/3)^{-p} N_1
\end{align*}
so we take $C_{1,n} :=  (\ve n/3)^{-p}N_1$, which clearly satisfies the summability condition.

For estimating $\P[A_{2,n}]$, we proceed similarly:
\begin{align*}
    \P[A_{2,n}] &\leq \sum_{1 \leq n_1 < n_2 \leq n} \P[ |\vp_{n_1}| > n^\gamma\text{ and } |\vp_{n-2}|> n^\gamma]\\
    &\leq \sum_{1 \leq n_1 < n_2 \leq n} \E[\P[ |\vp_{n_1}| > n^\gamma\text{ and } |\vp_{n-2}|> n^\gamma\,|\,\mathcal{F}_{n_2-1}]]\\
    &\leq \sum_{1 \leq n_1 < n_2 \leq n} \P[ |\vp_{n_1}| > n^\gamma]\P[|\vp_{n_2}| > n^\gamma]\\
    &\leq n^{2-2\gamma p}N_1^2
\end{align*}
Hence, we take $C_{2,n}:= n^{2-2\gamma p}N_1^2$, which also clearly satisfies the summability condition.

Note that because $\vp_n$ is a martingale difference $\E[\E[\vp_{n,k}\,|\,\mathcal{F}_{k-1}]] = -\E[\E[\vp_k - \vp_{n,k}\,|\,\mathcal{F}_{k-1}]]$ and so in particular these have the same magnitude.
\begin{align*}
    |T_n - \overline{T_n}| &= \Big|\sum_{1 \leq k \leq n}\E[\vp_{n,k}\,|\,\mathcal{F}_{k-1}]\Big|\\
    &\leq \sum_{1 \leq k \leq n} |\E[\vp_k - \vp_{n,k}\,|\,\mathcal{F}_{k-1}]|\\
    &\leq \sum_{1 \leq k \leq n} |\E[\vp_k1_{\{\vp_k \geq n^\gamma\}}\,|\,\mathcal{F}_{k-1}]|
\end{align*}
In particular, if $n$ is such that $\E[\vp_k1_{\{\vp_k \geq n^\gamma\}}\,|\,\mathcal{F}_{k-1}] \leq \ve/3$ for all $k < n$, $\P[A_{3,n}] = 0$; by definition of $N_2$ we immediately obtain $\P[A_{3,n}] \leq 1_{\{n \leq N_2\}}$, which satisfies the summability condition.

Finally, to bound $\P[A_{4,n}]$, first note that $\overline{\vp_{k,n}}$ is (for fixed $n$) a martingale difference and hence $\overline{T}_n$ is a martingale. This allows us to use the Burkholder's inequality relating the expectation quadratic variation to the expectation of $\overline{T}_n$. Specifically, we let $\overline{Q}_n = \sum_{1 \leq k \leq n} \overline{\vp_{n,k}}$. We fix $M$ to be the smallest integer at least $\frac{p}{2(1-\gamma)}$. There is, by the Burkholder's inequality, a constant $C$ (depending on $M$ but otherwise uniform) such that
\[
    \frac{1}{C} \E[\overline{Q}_n^M] \leq \E[\overline{T}_n^{2M}] \leq C \E[\overline{Q}_n^M] 
\]
yielding via Chebyshev
\begin{equation}
    \P[A_{4,n}] \leq Cn^{-2M} \E[\overline{Q}_n^M]
\end{equation}
so that in particular it now remains only to estimate the $M$-th moment of the quadratic variation.

We can write $\E[\overline{Q}_n^M]$ can be written as a linear combination of terms of the form $\E[\vp_{k_1,n}^{2q_1}\cdots \vp_{k_\ell}^{2q_\ell}]$, where $k_1<\cdots<k_\ell$ and $\sum_{k=1}^\ell q_k =  M$. We now fix $p_0 = \min\{p,2\}$. Recall that $|\overline{\vp}_{n,k}| \leq n^\gamma$, which gives us the bound
\[ \overline{\vp}_{n,k}^{2q} \leq n^{\gamma(2q -p_0)}|\overline{\vp}_{n,k}|^{p_0}\]
Because the process is a martingale, we exploit conditional independence and the bound $\E[|\overline{\vp}_{n,k}|^{p_0}] \leq N_1$ to obtain:
\begin{equation}
    \E[\vp_{k_1,n}^{2q_\ell}\cdots \vp_{k_m}^{2q_\ell}] \leq n^{\gamma(2M-\ell p_0)} N_1^\ell
\end{equation}
By well known estimates on multinomial coefficients and basic counting arguments, one sees that the number of terms corresponding to any given $\ell$ dominated by $Cn^\ell$, where $C$ is a constant which depends on $M$. Summing over all $\ell \leq M$, one obtains
\[ \E[\overline{Q}_n^M] \leq C \max\{N_1, N_1^M\}n^{2M\gamma}\]
where again, $C$ is a constant depending on $M$. Finally, by Chebyshev,
\[\P[A_{4,n}] \leq C n^{-2(1-\gamma)M}\]
for a constant $C$ depending on $M$ and on $N_1$.
\end{proof}

One sees easily that (\ref{phidom}) implies (\ref{unimomentcontrol}) and (\ref{unitailcontrol}).

The Wasserstein topologies we've introduced topologize two important quantities; closeness in this topology means closeness of the moments and closeness of the tails. Thus from our explicit version of Benoist and Quint's variant of the Baum--Katz theorem for martingales, we obtain the following:
\begin{thm}\label{realmartingalethm}
	Fix $p>1$, and let $K \subset \mathcal{P}^{p}(\R)$ be compact in the $W^p$ topology. For any $\ve > 0$, there are constants $C_n=C_n(\ve,K)$ satisfying $\sum_{n\geq 1}n^{p-2}C_n < \infty$ such that
	\[\P[|S_n|>\ve n] \leq C_n\]
	for the sum process associated to any martingale difference $\vp_n$ with law $\mathcal{L}(\E[\vp_n\,|\,\mathcal{F}_{n-1}]) \in K$ for all $n$.
\end{thm}

We postpone the proof of \Cref{realmartingalethm} until \Cref{invsec}, since it relies on ideas introduced in the proof of \Cref{logcompact2unif}.

Via this theorem it is possible to more or less immediately obtain the necessary estimates for our specific results on random walks on linear groups.

Estimates for martingales with a subexponential moment of the type we need are essentially proven in \cite{FAN2017538}; the following is a straightforward consequence of \cite[Theorem 2.1]{FAN2017538}. (See equation (1.5) in \textit{loc. cit.} and the surrounding discussion.)
\begin{prop}\label{semilogmgbd}
	Let $(\vp_n)_{n\geq 1}$ be a martingale difference and $S_n$ as defined in (\ref{summ}), such that
	\begin{equation}\label{subextail} \P[|\vp_n > t \,|\, \mathcal{F}_n] \leq \alpha e^{-t^\delta}\end{equation}
	for some $\alpha > 0$ and $0 < \delta < 1$. Then there exist positive $C$ and $c$ (depending only on $\alpha, \delta$ and $\ve$) such that
	\begin{equation}\P[|S_n| > n\ve] \leq C e^{-cn^{.99\delta}}\end{equation}
\end{prop}

In particular, in the subexponential moment regime, things are less delicate than in the fractional moment regime and control of the moments suffices. The appearance of the $.99$ above is a consequence of the necessity of a bound of the form

\[ \sup_n \E[\vp_n^2 \exp((\vp_n^+)^{\kappa})\,|\,\mathcal{F}_{n-1}] < \infty,\]
where $\vp_n^+:=\max\{\vp_n,0\}$; in particular the bound (\ref{subextail}) implies such bounds for all $\delta' < \delta$.

The following follows by standard Chernoff bounds:
\begin{prop}\label{expmgbd}
    Let $(\vp_n)_{n\geq 1}$ be a martingale difference and $S_n$ as defined in (\ref{summ}), such that
	\begin{equation}\label{subextailhyp} \P[|\vp_n > t \,|\, \mathcal{F}_n] \leq \alpha e^{-\beta t}\end{equation}
	for some $\alpha,\beta  > 0$. Then there exist positive $C$ and $c$ (depending only on $\alpha, \delta$ and $\ve$) such that
	\begin{equation}\P[|S_n| > n\ve] \leq C e^{-cn}\end{equation}
\end{prop}

\subsection{Uniform estimates for Markovian random walks on $G$-spaces}\label{bllnsec}
The work in the previous section is more or less all that is needed to control martingales which come from distributions in small neighborhoods in our topologies, so long as our martingale decompositions are well-behaved. We now show that this is the case, and moreover that the coboundary term and the regularity of the invariant measure both satisfy certain uniform estimates. Towards that end, we need uniform large deviation estimates for the Breiman law of large numbers as it applies to the random walks under consideration. The necessary result (and its proof) are very similar to \cite[Proposition 3.1]{bqclt}. Ensuring uniformity requires that we state the proposition in a more narrow context.

Given a topological group $G$ and a compact metrizable $G$-space $X$, i.e. a space equipped with a continuous $G$-action, any measure on $G$ naturally induces a Markovian random walk on $X$. Indeed, we just take for the kernel $P(x,A)=\mu[\{g \in G\,|\,gx \in A\}]$. We will also allow $P$ denote the operator on $C(X)$ acting by $Pf(x) = \int f(gx) \,d\mu(g)$.

Given any choice of initial distribution $\nu \in \mathcal{P}(X)$, we get a probability distribution on the space of forward trajectories $X^{\N}$, which we denote $\P_\nu$. We let $\overline{x} = (x_n)_{n\in\N}$ denote elements of $X^{\N}$. If $\nu = \delta_x$ for some $x$, we denote it by $\P_x$ rather than $\P_{\delta_x}$. The following was shown (in greater generality than we formulate here) in \cite{breiman}:

\begin{fact}[\cite{breiman}]
If $f$ has a unique average with respect to $\mu$, then ($\P_x$ almost surely, for any $x$) we have
\[ \lim_{n \rightarrow \infty} \frac{1}{n}\sum_{k=1}^n f(x_k) \rightarrow  \int f \,d\nu\]
(where $\nu \in \mathrm{Inv}(\mu)$).
\end{fact}

A key ingredient in the proof of a central limit theorem for distributions under a second moment assumption in \cite{bqclt} was the existence of large deviation estimates for this law. Their precise result is more general than what we formulate here. As we mentioned previously, their result doesn't require a unique average, whereas we were not able to uniformize this result, and moreover their result holds for general Markov--Feller chains rather than those induced by this kind of random walk on a group. See \cite[Proposition 3.1]{bqclt} for the precise statement; we avoid stating it here in order not to have to introduce all the general notions associated to Markov--Feller chains and their associated transfer operators.

\begin{lem}\label{unifblln}
    Let $G$ be a topological group, and $X$ a compact metrizable $G$-space. Let $\P_{x,\mu}$ denote the probability on the space of trajectories starting at $x$ coming from the Markovian random walk on $X$ induced by $\mu$. Fix $f \in C(X)$. If $K \subset \mathcal{P}(G)$ compact in the weak topology and
    \begin{enumerate}
        \item $K$ is tight; i.e. for any $\ve$ there is some compact set $K_\ve \subset G$ such that $\mu[K_\ve] \geq 1-\ve$ for all $\mu \in K$
        \item $f$ has a unique average for all $\mu \in K$
    \end{enumerate}
    Then for any $\ve > 0$, there are positive constants $C= C(K,f,\ve)$ and $c(K,f,\ve)$ such that
    \begin{equation}
        \P_{x,\mu}\left[ \left|\frac{1}{n}\sum_{k=1}^n f(x_n) - \int f \,d\nu_\mu\right| > \ve\right] \leq Ce^{-cn}
    \end{equation}
    for any $x \in X$, $\mu \in K$.
\end{lem}

Before undertaking the proof proper, we need to introduce certain notions and prove a technical lemma.
\begin{deffo}\label{almostinvMeasures}
Let $G$ be a topological group and $X$ a compact metrizable $G$-space. Given a measure $\mu$ on $G$ and $\nu$ on $X$, we say $\nu$ is \textit{$(\ve,\mu)$-almost invariant} if for any $f \in C(X)$ we have
\begin{equation}
    \int f d\nu - \int f d(\mu\ast \nu) \leq \ve \|f\|_\infty
\end{equation}
\end{deffo}
The necessary technical lemma is the following:
\begin{lem}\label{uniformuniformconv}
    Let $G$ be a topological group, and $X$ a compact metrizable $G$ space. Fix $f \in C(X)$. Let $K \subset \mathcal{P}(G)$ be a compact set in the weak topology  such that
    \begin{enumerate}
        \item The family $K$ is tight,
        \item $f$ has a unique average for all $\mu \in K$ (see \Cref{def:uniqueAvg}).
    \end{enumerate}Then there is $C = C(K,f)$ such that for any $(\ve, \mu)$-almost invariant measure $\nu'$ on $X$ with $\mu \in K$
        \begin{equation}
            \left|\int f d\nu' - \int fd\nu_\mu \right| < C\ve
        \end{equation}
        where $\nu_\mu$ is an invariant measure associated to $\mu$.
\end{lem}
\begin{proof}
    Throughout, for some $\mu \in K$, we let $\nu_\mu$ denote an invariant measure associated to it. Note that by our assumption of a unique average, $\int f\,d\mu_\nu$ does not depend on the arbitrary choice of invariant measure. We argue by contradiction, assuming the existence of $(\mu_n,\nu_n) \in K \times \mathcal{P}(X)$ such that the $\nu_n$ are $(2^{-n}, \mu_n)$-almost invariant, and 
    \[\int f d\nu_n \geq  \int fd\nu_{\mu_n} + 1\]
    (Either such a sequence exists, or there is a symmetric case which one treats identically.)
    
    By compactness, we can pass to a subsequence converging to some $(\mu,\nu)$. Clearly
    \begin{align*}
        \int f\,d\nu 
        &\geq \limsup \int f\,d\nu_{\mu_n} + 1\\
        &\geq \int f d\nu_{\mu} + 1
    \end{align*}
    The last line follows by the fact that the invariant measures of the $\mu_n$ converge to an invariant measure for $\mu$. Clearly this cannot hold if $\nu$ is invariant (i.e. $\nu = \nu_\mu$), so it suffices to demonstrate said invariance. This in turn follows easily from showing the weak convergence of $\mu_n\times \nu_n$ to $\mu \times \nu$.

    Note that because the sequence $(\mu_n,\nu_n)$ is tight, the desired weak convergence is equivalent to the following:
    \begin{equation}
        \int h \,d(\mu_n \times \nu_n)\rightarrow \int h \,d(\mu \times \nu)
    \end{equation}
    for all $h \in C_0(G\times X)$, the space of all continuous functions on $G\times X$ vanishing at infinity. Importantly, we don't need to consider the entire space of bounded continuous functions. In particular, in the space of functions vanishing at infinity, the algebra generated by functions of the form $h(g,x) = h_1(g)h_2(x)$ for $h_1 \in C_0(G)$, $h_2\in C(X)$ is dense.

    Hence for any $h \in C_0(G\times X)$, we can find a decomposition $h = h_1+ h_2$ where $h_1$ is a linear combination of products of functions of one variable, and $\|h_2\|_\infty < \ve$.  Clearly 
    \[\int h_1\,d(\mu_n\times \nu_n) \rightarrow \int h_1\,d(\mu\times \nu)\]
    and so in particular because all measures involved are probability measures,
    \[ \limsup_{n\rightarrow \infty} \int h\,d(\mu_n \times \nu_n) - \int h \,d(\mu\times \nu) < 2\ve\]
    for any $\ve > 0$, showing the desired weak convergence, and hence $\mu$-invariance of $\nu$. This gives the desired contradiction, establishing the lemma.
\end{proof}

This has the following corollary, which we use in proving Lemma \ref{unifblln}.
\begin{cor}\label{unifunifconvcor}
    Let $G$, $K\subset \mathcal{P}(G)$, $X$ and $f$ satisfy the hypotheses of Lemma \ref{uniformuniformconv}. Then 
    \[
        \sup_{\mu,x} \left|\frac{1}{n}\sum_{k=1}^n \int f(gx) \,d\mu^{\ast k}(g) - \int f(x') d\nu_\mu(x')\right| \rightarrow 0\]
\end{cor}
Indeed, it follows immediately from the fact that the measure $\frac{1}{n}\sum_{k=1}^n \mu^{\ast k}\ast \delta_x$ is $(\frac{2}{n},\mu)$-almost invariant.

\begin{proof}[Proof of \Cref{unifblln}]
    Fix $K, f$ and $\ve$ as set out in  the hypotheses of the theorem. Recall that we let $Pf(x) = \int f(gx)\,d\mu(g)$. By Corollary \ref{unifunifconvcor}, there is $N = N(K,f,\ve)$ (and not depending on $x, \mu \in K$) such that $n>N$ implies

    \begin{equation}
        \left|\frac{1}{n}\sum_{k=1}^n P^kf(x) - \int f(x') d\nu_\mu(x')\right| < \ve
    \end{equation}
    for any $x \in X$, $\mu \in K$. From here, we could apply the original strategy of Benoist and Quint; to keep this paper complete without reproducing their work, we instead appeal to work of Cai, Duarte and Klein \cite[Theorem 1.1]{cai2024statisticalpropertiesmixingmarkov} to immediately extract the desired large deviation estimates. (We have technically not shown that the hypotheses of this result are fulfilled, but by \cite[Remark 1.2]{cai2024statisticalpropertiesmixingmarkov} we can nevertheless conclude the existence of our desired large deviation estimates.)
\end{proof}

\begin{remm}
    The results of Cai, Duarte, and Klein used here are effective in the sense that given effective mixing estimates, one obtains effective constants. Our results are not effective because our input into their machinery is not effective. In particular, we use a compactness argument. Using e.g. spectral methods, it is well known that for sufficiently regular observables, one can get effective uniform mixing estimates (the necessary input) if one has effective bounds on the rate of convergence of the Birkhoff averages; this relies on the spectral theory of the transfer operator associated to the Markov chain, restricted to certain classes of functions.

    It seems plausible that one should be able to obtain effective estimates of this kind at the very least in the semi-exponential regime and for all local fields; using \cite[Theorem 1.1]{cai2024statisticalpropertiesmixingmarkov} would then yield effective large deviation estimates in such regimes. In particular, if these estimates only depend on $\mu$ via its moments, then one should be able to replace the compactness assumption with a boundedness assumption in our main abstract results. Depending on the assumptions made, it is possible that such results could even give uniform large deviation estimates which hold even outside the locus of $\mu$ which have unique invariant measures/unique averages.
\end{remm}

Note that for the case of $G = \mathrm{SL}_d(\mathbb{K})$ with $\mathbb{K}$ a local field, compactness in any of the Wasserstein type topology studied in this work implies $K$ satisfies the weak-compactness and tightness hypotheses --- the Heine--Borel property holds for these groups equipped with the appropriate metric and so by Chebyshev, one obtains tightness. This suffices for our large deviation estimates, but is in fact the case that any open $\ve$ ball in a Wasserstein topology is weakly precompact; this is crucial for our proof of continuity of the variance appearing in the central limit theorem. (This is an immediate consequence of the aforementioned tightness and Prokhorov's theorem.)

\section{Large deviations for linear cocycles and estimates on the invariant measure}\label{invsec}
From the large deviation estimates obtained for cocycles and for martingales, one immediately obtains parts of Theorems \ref{logmatrixlde} and \ref{semilogmatrixlde}, specifically the results concerning $\log \|A_n^\mu x\|$ by taking $G=\GL_d(\KK)$, $X = \Proj(\mathbb{K}^d)$, and $\Phi(M,x) = \log\frac{\|Mx\|}{\|x\|}$. Recall that throughout we will be abusing notation, identifying a direction $x \in \Proj(\mathbb{K}^d)$ with an arbitrary vector in said direction whenever no confusion is likely.

We require a certain technical lemma:
\begin{lem}\label{logcompact2unif}
    Given $K \subset \mathcal{P}_{\log}^p$ compact for some $p >1$, the following hold:
    \begin{enumerate}
        \item For any $x \in \Proj(\mathbb{K}^d)$, we have \begin{equation}\label{uniflogmomentseven} \sup_{\mu \in K} \int \left|\log\frac{\|Mx\|}{\|x\|}\right|^p \,d\mu(M) < \infty \end{equation}
        \item For any $x \in \Proj(\mathbb{K}^d)$ and $\ve > 0$, there is $N_2 = N_2(\ve,K)$ such that
        \begin{equation}
            \int 1_{\left\{\left|\log \frac{\|Mx\|}{\|x\|}\right| \geq N_2\right\}}\left|\log\frac{\|Mx\|}{\|x\|}\right|\,d\mu(M) < \ve 
        \end{equation}
    \end{enumerate}
\end{lem}

\begin{proof}
    Clearly it suffices to prove this holds in a neighborhood of any point in $\mathcal{P}^p_{\log}$; the fact that the topology is metrizable means in particular we can work sequentially. It suffices to show that if $\mu_n \overset{W^p_{\log}}{\rightarrow}\mu$, these properties holds for $\mu_n$ when $n$ is sufficiently large.

    We proceed along these lines; the moment bound follows more or less immediately by Proposition \ref{momentmetric}, as the logarithmic moments converge of $\mu_n$ converge to those of $\mu$ and in particular are eventually less than twice that of $\mu$.

    To establish the second condition, tail control, we exploit the fact that convergence in any of the concave Wasserstein topologies implies convergence in the weak topology.

    To simplify notation, we fix $x \in \Proj(\mathbb{K}^d)$ and then abusively identify $\mu$ with its pushforward onto $\R$ under the map $M \mapsto \left|\log \frac{\|Mx\|}{\|x\|}\right|$. Specifically, we find $N_2 >0$ such that $\int_{t\geq N_2} t^p d\mu(t) < \ve/12$. Then we define $f$ by
    \[ f(t) = \begin{cases} t^p \quad&\text{for}\quad t \leq N_2\\
    N_2^p-N_2^p(t-N_2)\quad&\text{for}\quad t \in (N_2, N_2+1)\\
    0\quad&\text{for}\quad t \geq  N_2+1\end{cases}\]
    We let $I_1 = t^p1_{\{ t > N_2\}}$ and $I_2 = t^p1_{\{t > N_2+1\}}$. In particular, $0 \leq I_2 \leq f \leq I_1$ and $f$ is continuous. It suffices to show that $\int I_2 \,d\mu_n < \ve$ for $n$ sufficiently large. We split the integral as follows:
    \begin{align*}
        \int I_2\,d\mu_n &= \int f \, d\mu_n + \int I_2 - f \,d\mu + \int I_2 - f\,d(\mu_n-\mu).
    \end{align*}
    The first term is small by assumption on $N_2$; the second term goes to zero by dominated convergence. For the last term, note that $|\int (I_2-f)\,d(\mu_n - \mu)| \leq |\int I_2 - I_1\,d(\mu_n-
    \mu)|$. Finally, since $I_2 - I_1$ is continuous when restricted to the support of the function, by the weak convergence this term also decays to zero.
\end{proof}
This is enough to prove Theorems  \ref{logmatrixlde}, \ref{semilogmatrixlde}, and \ref{unifldeexp}, but we will prove all our large deviation theorems at once. We will take for granted these theorems for now, and prove certain estimates on the invariant measure, which are necessary for the proof of \Cref{matrixcoefflde}. We can now give the proof of \Cref{realmartingalethm}:

\begin{proof}[Proof of \Cref{realmartingalethm}]
	It suffices to show that for any $\mu \in \mathcal{P}^{p}(\R)$, there is a corresponding $W^p$ neighborhood $\mathcal{N}$ such that one can find, for all $\ve > 0$, constants $\tilde{C}_n:= \tilde{C}_n(\ve,N)$ which satisfy $\sum_{n\geq 1}n^{p-2}\tilde{C}_n$ and
	\[ \P[|S_n|>n\ve]\]
	for $\vp_k$ a martingale difference with law $\mathcal{L}(\E[\vp_n\,|\,\mathcal{F}_{n-1}]) \in \mathcal{N}$. Towards this end, we need only show that on sufficiently small $W^p$ neighborhood, (\ref{unimomentcontrol}) and (\ref{unitailcontrol}) hold. This essentially follows by the argument laid out in \Cref{logcompact2unif}. 
\end{proof}

\subsection{Estimates on the invariant measure}
The proof of the corresponding fact in the semi-logarithmic regime follows by the same argument, or more precisely by the first half of the same proof, tail properties not mattering in this regime, and the following clear fact:
\begin{lem}\label{semilogcompact2unif}
    Given $K \subset \mathcal{P}_{\slog}^\delta(\SL_d(\KK))$ compact for some $\delta \in (0,1)$, then
    \begin{equation}\label{unifsemilogmoment} \sup_{\mu \in K} \int \exp\left(\left|\log\frac{\|Mx\|}{\|x\|}\right|^\delta \right)\,d\mu(M) < \infty \end{equation}
\end{lem}

\begin{remm}
    Compactness suffices for our applications, but note that in the semi-logarithmic regime (i.e. in Lemma \ref{semilogcompact2unif}), it in fact suffices to take $K$ a bounded set in the metric, at least as far as the martingale part is concerned. Our method as it stands cannot show the same for the part concerning Markov chains, and so our results currently still depend on a compactness assumption; we expect that this is more an artifact of the proof than a fundamental limitation in the semi-logarithmic regime. On the other hand, in the logarithmic regime, we needed to work over compacta to get the results we obtain (see \cite[Remark 3.2]{bqclt}).
\end{remm}

The following propositions are uniform versions of \cite[Proposition 4.1]{bqclt} and \cite[Corollary 4.2]{bqclt}; they essentially follow by combining our uniform results in Theorems \ref{logmatrixlde} and \ref{semilogmatrixlde} with the proofs of said results. The proofs are virtually identical to those presented by Benoist and Quint, and the only new inputs are Theorems \ref{logmatrixlde} and \ref{semilogmatrixlde}.

\begin{prop}\label{furstldesprop}
	Let $K \subset \mathcal{P}^{p}_{\log}$ for $p> 1$ (or $K \subset \mathcal{P}^\delta_{\slog}$ for $\delta \in (0,1)$) be compact, such that all the $\mu \in K$ are  irreducible.
	
	In the logarithmic regime, there are, for all $\ve > 0$, constants $C_n = C_n(K,\ve)$ satisfying $\sum_{n\geq 1} n^{p-2}C_n < \infty$ and
	\begin{equation}\label{logfurstnormlde}
		\P\left[|\log\|A_n^\mu\| - n\lambda_1(\mu)|>n\ve\right] \leq C_n.
	\end{equation}
	Moreover, for any non-zero $x \in \mathbb{K}^d$
	\begin{equation}\label{logfurstvectnormlde}
		\P\left[\left|\log\frac{\|A_n^\mu x\|}{\|x\|} - n\lambda_1(\mu)\right|>n\ve\right] \leq C_n.
	\end{equation}
	In the semi-logarithmic regime, one has $C = C(K,\ve)$ and $c=c(K,\ve)$ such that
	\begin{equation}\label{semilogfurstnormlde}
		\P\left[|\log\|A_n^\mu\| - n\lambda_1(\mu)|>n\ve\right] \leq Ce^{-cn^{.99\delta}}
	\end{equation}
	and moreover for any non-zero $x \in \mathbb{K}^d$,
	\begin{equation}\label{semilogvectnormlde}
		\P\left[\left|\log\frac{\|A_n^\mu x\|}{\|x\|} - n\lambda_1(\mu)\right|>n\ve\right] \leq Ce^{-cn^{.99\delta}}.
	\end{equation}
\end{prop}

\begin{prop}\label{lyaptwoprop}
	In the logarithmic regime, there are for all $\ve > 0$ constants $C_n = C_n(K,\ve)$ such that $\sum_{n\geq 1} n^{p-2}C_n <\infty$ and
	\begin{equation}\label{loglyaptwolde}
		\P[|\log\|\wedge^2 A_n^\mu\| - n(\lambda_1(\mu)+\lambda_2(\mu))|>n\ve] \leq C_n.
	\end{equation}
	In the semi-logarithmic regime, there are positive $C = C(K,\ve)$ and $c=c(K,\ve)$ such that
	\begin{equation}\label{semiloglyaptwolde}
		\P[|\log\|\wedge^2 A_n^\mu\| - n(\lambda_1(\mu)+\lambda_2(\mu))|>n\ve] \leq Ce^{-cn^{.99\delta}}.
	\end{equation}
\end{prop}

\begin{remm}
	If one drops the irreducibility condition for $\mu \in K$, the large deviation estimates for matrix norms (\ref{logfurstnormlde}) and (\ref{semilogfurstnormlde}) and are still valid in Proposition \ref{furstldesprop}; in particular one then obtains Proposition \ref{lyaptwoprop} as a corollary.
\end{remm}

We will also need the following:

\begin{prop}
	Let $K$ be compact in either $\mathcal{P}^p_{\log}(\SL_d(\mathbb{K}))$ (for $p > 1$) or $\mathcal{P}^{\delta}_{\slog}(\SL_d(\mathbb{K}))$ (for $\delta \in (0,1)$), such that all $\mu \in K$ are strongly irreducible and contracting. Then \begin{equation}\inf_{\mu \in K} \lambda_1(\mu) > 0.\end{equation} If $\wedge^2_\ast (K)$ is compact in the topology induced by $W^1_{\log}$ we also obtain \begin{equation}\inf_{\mu \in K} \lambda_1(\mu)-\lambda_2(\mu) > 0. \end{equation}
\end{prop}

These two bounds follow almost immediately by existing work; by work of Furstenberg, $\lambda_1(\mu) > 0$ for all $\mu \in K$, and by work of Furstenberg and Kifer, $\lambda_1(\mu)$ is continuous in $\mu$. By work of Guivarc'h, $\lambda_1(\mu)-\lambda_2(\mu) > 0$ for all $\mu \in K$, and by Theorem \ref{usclyap}, said quantity is lower semi-continuous. Recall that we have shown any $K$ compact in a semi-logarithmic Wasserstein topology has compact image in $W^1_{\log}$ under the map $\mu \mapsto \wedge^2_\ast \mu$ in Theorem \ref{wedgecontinuity}, so that in particular this condition is trivial in said context. (Recall also that this last fact is the content of \Cref{unifgapbound}.)

We require a purely deterministic technical lemma, which requires the introduction of certain terminology. Our notation and terminology differ, but what follows is essentially \cite[Lemma 4.7]{bqclt}.

For any $M \in \SL_d(\mathbb{K})$, $M$ admits a ``KAK'' factorization:
\begin{equation}
	M = K_MA_M\tilde{K}_M
\end{equation}
where $K_M, \tilde{K}_M$ are isometries, and $A_M$ is diagonal, of the form $A_M = \mathrm{diag}(a_1,\cdots,a_d)$ where $|a_1| \geq |a_2| \geq \cdots \geq |a_m|$. We let $\omega(M) \in \Proj(\mathbb{K}^d)$ denote the ``outgoing density point'' of $M$, defined by
\begin{equation}
	\omega_M := \mathbb{K}K_Me_1
\end{equation} 
Similarly, we let $\iota(M) \in \Proj((\mathbb{K}^d)^\ast)$ denote the ``incoming density point'' of $M$, defined by
\begin{equation}
	\iota_M := \mathbb{K}^\ast e_1^\ast \circ \tilde{K}_M
\end{equation}
The significance of these density points is as follows: for a given $M$, $\iota(M)$ is essentially a projection onto the ``incoming'' direction of most increase (essentially because of the projectivization), and $\omega(M)$ is the ``outgoing'' direction of most increase. If $|a_1| > |a_2|$, there is a unique direction of fastest increase; this will be the case asymptotically for the matrices under consideration due to our dynamical assumptions.

For an element $x \in \Proj(\mathbb{K}^d)$ and $f \in \Proj((\mathbb{K}^d)^\ast)$, we let \[|f(v)|= \frac{|f(v)|}{\|f\|\|v\|}\]
where we abuse notation to let $f$ and $v$ denote simultaneously elements of projectivized space and non-zero elements of the corresponding equivalence classes. 

We also define a distance on projective space $\Proj(\mathbb{K}^d)$ by:
\[ d(x_1,x_2) = 
\frac{\|x_1\wedge x_2\|}{\|x_1\|\|x_2\|}\] where we equip $\bigwedge^2 \KK^d$ with e.g. the projective norm detailed in e.g. \cite[Appendix 4]{quas2019explicit}. (If one equips the space of finite formal sums of pairs $(x_1,x_2)$ with the norm induced by $\| \sum c_i(x_{1,i},x_{2,i} )\| = \sum c_i \|x_1\| \cdot \|x_2\|$, then this projective norm is the quotient norm associated to taking the induced norm structure on the tensor product $\bigotimes^2 \KK^d$ and then the quotient necessary to obtain $\bigwedge^2 \KK^d$.)
\begin{remm}
    In the case where $\KK = \R$ or $\C$, if one has an inner product structure inducing the norm on $\KK^d$ then one can take a norm on $\bigwedge^2 \KK^d$ such that
\[ d(x_1,x_2)^2 = 1- |\langle x_1,x_2\rangle|^2\]
and in the Archimedean setting one can work with this norm instead of the projective norm; all results below still hold. This is precisely the metric considered in e.g. \Cref{regularitythmintro}.
\end{remm}
Finally, for any $M$ we define
\[ \gamma(M) = \frac{\|\wedge^2 M\|}{\|M\|^2} = \frac{|a_2|}{|a_1|} \]
which we will call the \emph{multiplicative gap}. We have a unique direction of most increase precisely when $\gamma(M) < 1$. We can now formulate the following, which is proven by Benoist and Quint: 
\begin{lem}[{\cite[Lemma 4.7]{bqclt}}]
	For any $x \in \Proj(\mathbb{K}^d)$, $f \in \Proj((\mathbb{K}^d)^\ast)$ and $M \in \SL_d(\mathbb{K})$ we have the following estimates:
	\begin{equation}\label{x2ingobd}\iota_M(x) \leq \frac{\|Mx\|}{\|M\|\|x\|} \leq \iota_M(x) + \gamma(M)\end{equation}
	\begin{equation}\label{f2outgobd}  f(\omega_M) \leq \frac{\|f\circ M\|}{\|f\|\|M\|} \leq f(\omega_M)+ \gamma(M)\end{equation}
	\begin{equation}\label{timesbd} d(Mx, \omega(M)) \iota_M(x) \leq \gamma(M)
	\end{equation}
	(Recall that we abuse notation, identifying elements of projective space with non-zero vectors in the corresponding equivalence classes throughout; all relevant quantities are invariant under non-zero scaling.)
\end{lem}

The essential technical results for the regularity result in the semi-logarithmic regime and the large deviation estimates for matrix coefficients in both regimes are probabilistic variants of these three estimates. In the logarithmic moment regime, which is essentially a uniform version of \cite[Lemma 4.8]{bqclt}:

\begin{lem}\label{logalignbd}
	Let $K \subset \mathcal{P}_{\log}^p(\SL_d(\mathbb{K}))$ be compact such that all $\mu \in K$ are strongly irreducible and contracting. Then for any $\ve > 0$, there are constants $C_n = C_n(K,\ve)$ satisfying $\sum_{n\geq 1} n^{p-2}C_n < \infty$ such that for all $\mu \in K$, $x \in \Proj(\mathbb{K}^d)$ and $f \in \Proj((\mathbb{K}^d)^\ast)$ we have
	\begin{equation}
		\P[d(A_n^\mu x,\omega(A_n^\mu)) \geq e^{(\lambda_1(\mu)-\ve )n}] \leq C_n
	\end{equation}
	\begin{equation}
	\P[ f(\omega(A_n^\mu))\leq e^{-\ve n}] \leq C_n 
	\end{equation}
	\begin{equation}\label{logbigalign}
	\P[f(A_n^\mu x) \leq e^{-\ve n}] \leq C_n 
	\end{equation}

\end{lem}

In the semi-logarithmic regime, we have the same but with stronger estimates:

\begin{lem} \label{semilogalignbounds}
	Let $K \subset \mathcal{P}_{\log}^p(\SL_d(\mathbb{K}))$ be compact (with compact image in the $W^1_{\log}$ topology under the map $\wedge^2_\ast$) such that all $\mu \in K$ are strongly irreducible and contracting. Then for any $\ve > 0$, there are positive $C=C(K,\ve)$ and $c=c(K,\ve)$ such that for all $\mu \in K$, $x \in \Proj(\mathbb{K}^d)$ and $f \in \Proj((\mathbb{K}^d)^\ast)$ we have
	\begin{equation}
		\P[d(A_n^\mu x,\omega(A_n^\mu)) \geq e^{-(\lambda_1(\mu)+\ve) n}] \leq Ce^{-cn^{.99\delta}}
	\end{equation}
	\begin{equation}
		\P[\langle f, \omega(A_n^\mu)\rangle \leq e^{-\ve n}] \leq Ce^{-cn^{.99\delta}}
	\end{equation}
	\begin{equation}\label{semilogbigalign}
		\P[f(A_n^\mu x) \leq e^{-\ve n}] \leq Ce^{-cn^{.99\delta}}
	\end{equation}
	
\end{lem}

The proofs of these two lemmas are virtually identical, and based off the proof of \cite[Lemma 4.8]{bqclt}; we will prove Lemma \ref{semilogalignbounds} in detail.
\begin{proof}
    By Theorem \ref{semilogmatrixlde}, for any choice of $\ve'> 0$, with probability at least $1-Ce^{-cn^\delta}$, we have the following four bounds, for any $x \in \Proj(\mathbb{K}^d)$ and $f \in \Proj((\mathbb{K}^d)^\ast)$
    \begin{align*}
        |\lambda_1(\mu)- \frac{1}{n}\log \|A_n^\mu\|| &\leq \ve'\\
        |\lambda_1(\mu)- \frac{1}{n}\log \|A_n^\mu x\|| &\leq \ve'\\
        |\lambda_1(\mu)- \frac{1}{n}\log \|f \circ A_n^\mu \|| &\leq \ve'\\
        |\lambda_1(\mu)+\lambda_2(\mu)- \frac{1}{n}\gamma_1(A_n^{\mu})| &\leq \ve'
\end{align*}
for any choice of $\ve'$ and all $\mu \in K$. It suffices to show that for arbitrarily large $n$ (uniformly in $\mu$, $x$, $f$) the above inequalities imply the sought bounds provided we choose $\ve'$ appropriately. Henceforth we presume the above bounds to hold. By (\ref{x2ingobd}), we obtain
\begin{align*}
\langle \iota(A_n^\mu), A_n^\mu x \rangle &\geq \frac{\|A_n^\mu x\|}{\|A_n^\mu\|\|x\|} - \gamma_1(A_n^\mu)\\ 
&\geq \exp(-2\ve' \gamma_{\min}n) - \exp((-1-\ve')\gamma_{\min}n)\\
&\geq \exp(-3\ve'\gamma_{\min}n)
\end{align*}
for sufficiently large $n$, requisite largeness depending only on $\lambda_{\min}$. Then, using (\ref{timesbd}), we readily obtain
\[ d(A_n^\mu x, \omega_{A_n^\mu}) \leq \exp((1-4\ve')\gamma_{\min}n)\]
proving the first of the sought estimates. The second estimate follows from the argument for the bound
\[\iota_{A_n^\mu}(x) \geq \exp(-3\ve'\gamma_{\min}n)\]
by duality. Finally, using the basic estimate
\[ f(x_1)\geq f(x_2) - d(x_1,x_2)\]
we obtain
\begin{align*}
     f(A_n^\mu x) &\geq f( \omega(A_n^\mu) ) - d(\omega(A_n^\mu), A_n^\mu x )\\
    &\geq \exp(-3\ve'\gamma_{\min}n) - \exp(-(1+4\ve')\gamma_{\min}n)\\
    &\geq \exp(-4\ve' \gamma_{\min}n)
\end{align*}
Thus, $\ve' < \frac{\ve}{4\gamma_{\min}}$ suffices, concluding the proof.
\end{proof}

\subsection{Proofs of large deviation theorems for the random walk}
Here we put together all the pieces to prove our main large deviation theorems in the logarithmic, semi-logarithmic, and fractional moment settings:

\begin{proof}[Proof of Theorems \ref{logmatrixlde}, \ref{semilogmatrixlde}, and \ref{unifldeexp}]

By \Cref{logcompact2unif}, for any $K$ compact in $\mathcal{P}_{\log}^p$ (resp. $\mathcal{P}_{\slog}^p$ and $\mathcal{P}^p$) we have precisely the uniform estimates needed to apply \Cref{explicitmgpolylde}, which gives us (in the logarithmic case) precisely the necessary large deviation estimates on the martingale $\Phi'(M,x)$ appearing in (\ref{mgmarkovdecomp}).

In the semi-logarithmic setting, the necessary large deviation estimates are a consequence of \Cref{semilogmgbd}; in the fractional/polynomial setting, such large deviations follow immediately from \Cref{expmgbd}.

Finally, in all cases we have uniform exponential large deviation bounds for the Markov chain term in decomposition (\ref{mgmarkovdecomp}) by \Cref{unifblln}.
\end{proof}

\Cref{matrixcoefflde} follows more or less immediately from the large deviation estimates already established for $\log \|A_n^\mu x\|$ together with formulae (\ref{logbigalign}) and (\ref{semilogbigalign}). 

\begin{proof}[Proof of \Cref{matrixcoefflde}] Fixing $f \in (\mathbb{K}^d)^\ast$ and $x \in \mathbb{K}^d$ in the associated affine spaces, both with norm one, we obtain:
\begin{equation} \log |f(A_n^\mu x)| = \log \|A_n^\mu x\| + \log \frac{|f( A_n^\mu x)|}{\|A_n^\mu x \|} \end{equation}

Clearly
\begin{align*}  
&\,\P[|\log |f(A_n^\mu x)| - n\lambda_1(\mu)|> n\ve] \leq\P[|\log \|A_n^\mu x\| - n\lambda_1(\mu)|>n\ve/2] + \P\left[\left|\log \frac{|f(A_n^\mu x)|}{\|A_n^\mu x \|}\right| > n\ve/2 \right]\end{align*}

Having already established the large deviation bounds mentioned in Theorems \ref{logmatrixlde} and \ref{semilogmatrixlde}, it suffices to bound the second term, which we have accomplished in \Cref{semilogalignbounds} and \Cref{logalignbd}. This covers the logarithmic and semi-logarithmic regimes; one could easily carry out the argument for the fractional regime as well, but this regularity is already known anyways, see \cite{bougerol2012products}.
\end{proof}

 \subsection{Regularity results}
In this section, we prove various regularity results, the proofs of all of which are essentially modelled off the proof of the following theorem of Benoist and Quint, which is one of their main technical results:

\begin{thm}[{\cite[Proposition 4.5]{bqclt}}] 
    Let $\mu$ is a strongly irreducible and contracting distribution on $\SL_d(\KK)$ with a finite logarithmic moment of order $p\geq 2$, and let $\nu$ denote the associated invariant measure on $\P(\KK^d)$. Then
    \[ \int_{\P(\KK^d)} |\log |f(x)|\, |^{p-1} \,d\nu(x) \]
    is finite for all $f \in (\KK^d)^\ast$ with unit norm, and moreover
    $f \mapsto \int_{\P(\KK^d)} |\log f(x) | \,d\nu(x)$ is continuous.
\end{thm}

Benoist and Quint derive this result from the bounds we have collected in our \Cref{logalignbd}. We briefly sketch their argument, since ours is a straightforward modification.

By \Cref{logbigalign}, we have in particular (for fixed $\ve > 0$ and $f \in \P((\KK^d)^\ast))$ with $\|f\|=1$:
\begin{align*}
\mu^{\ast n}\left(\left\{g\, : \, |f(gx)| \leq e^{-\ve n} \right\}\right) = \P[ |f(A_n^\mu x)| \leq e^{-\ve n}] \leq C_n
\end{align*}
for $C_n$ such that $\sum n^{p-2} C_n < \infty$. By $\mu$ invariance, the above inequality in fact yields
\[ \nu\left(\left\{x\,:\, \log| f(x)| < - \ve n \right\}\right) \leq C_n \]
Then one has:
\begin{align*}
    \int \left|\log |f(x)|\right|^{p-1} d\nu(x) &\leq \sum_{n \geq 1} \ve^{p-1}n^{p-1} [\nu(\{|f(x)| \leq e^{-\ve n}\}) - \nu(\{|f(x)| \leq e^{-\ve (n-1)}\})]\\
    &\leq \ve^{p-1} + \ve^{p-1}\sum_{n\geq 1} [(n+1)^{p-1}-n^{p-1}]C_n\\
    &\leq \ve^{p-1}+(p-1)2^p\ve^{p-1}\sum_{n\geq 1} n^{p-2} C_n < \infty
\end{align*}
establishing finitude. Continuity of the function $f \mapsto \int |\log |f(x)||^{p-1}\,d\nu(x)$ follows by noting that this function is a uniform limit of its (continuous) truncations $f \mapsto \int \min\{|\log |f(x)||^{p-1},(\ve n)^{p-1}\}\,d\nu(x)$. That this convergence is uniform follows from the fact that the constants $C_n$ don't depend on $f$; in particular the differences are uniformly bounded by tails of $(p-1)2^p \ve^{p-1} \sum_{n \geq 1} n^{p-2}C_n$.

While not explicitly stated in \cite{bqclt}, this implies log-H\"older regularity of the invariant measure $\nu$ very straightforwardly. Before introducing the definition of log-H\"older (in our specific context), we recall that if we endow $\R^d$ or $\C^d$ with their standard inner product structures, we can define a natural metric on their projectivizations by
\[d(x,y) = 1 - |\langle x, y \rangle|^2\]
where $x,y$ are identified with normalized vectors in the appropriate directions in the inner product space.
\begin{deffo}
We say a measure $\nu$ on $\P(\KK^d)$ (for $\KK = \R$ or $\C$) is \textit{log-H\"older} if for all sufficiently small $\ve> 0$ and all $x \in \P(\KK^d)$, $\nu(B_\ve(x)) \leq C|\log(\ve)|^{-\alpha}$ for some $C,\alpha > 0$.
\end{deffo}
Of course, the notion of log-H\"{o}lder measures makes sense on general metric spaces. It is sufficient, and in fact necessary by compactness of $\P(\KK^d)$, to exhibit a uniform bound on  
\[ \int_{\P(\KK^d)} |\log d(y,x)|^\alpha\,d\nu(x) \]
for all $y$. In particular, modulo the very simple technical \Cref{distanceslem} to follow, we obtain the following as a corollary.

\begin{cor}
    Given any $\mu$ on $\GL_d(\R)$ or $\GL_d(\C)$ with a finite logarithmic moment of degree $p \geq 1$ and moreover strongly irreducible and contracting, the invariant measure $\nu$ is log-H\"older.
\end{cor}

Indeed, it suffices to demonstrate that small $\ve$ balls $B_\ve(x)$ in the specified distance are contained in sets such that $f(x)$ is small for an appropriate functional:

\begin{lem}\label{distanceslem}
    For any $y$ there is a linear functional $f$ of unit norm such that $d(x,y) \geq \frac{|f(x)|^2}{2}$.
\end{lem}
\begin{proof}
    It suffices to take $f$ acting by $y \mapsto \langle z, \cdot \rangle$ for any $z$ orthogonal to $x$; since $d(z,x) = 1$, $\sqrt{1 - |\langle z, y \rangle|^2} = d(z,y) > 1-d(x,y)$. Using the simple bound
    \[ \sqrt{1-x} \leq 1-x/2\]
    valid for $x \in (0,1)$ gives the result.
\end{proof}
Then the log-H\"older property is an immediate consequence of the bounds on $\int |\log |f(x)||^{p-1}\,d\nu(x)$. 
\begin{remm}
    It was shown in \cite{monakov2024log} that the log-H\"older property in very high generality for random transformations acting on a Riemannian manifold; their approach also yields a better constant $\alpha$ in this special case when $p<2$; one can then obtain $\alpha = p/2$.
\end{remm}

Similarly, \Cref{regularitythmintro} is a straightforward corollary of the following theorem:
\begin{thm}\label{semilogfunctionalreg}
    Let $\mu$ be strongly irreducible and contracting, with a finite semi-logarithmic moment of order $\delta \in (0,1)$, and let $\nu$ denote the associated invariant measure on $\P(\KK^d)$. Then for any $\rho < .99\delta$,

    \begin{equation}\label{semilogint2bound}\int_{\P(\KK^d)} \exp\left(\left(\log \frac{1}{|f(x)|}\right)^\rho \right)\,d\nu(x)\end{equation}
    is uniformly bounded over all $f \in (\KK^d)^\ast$ of unit norm.
\end{thm}

\begin{proof}
    We essentially substitute the stronger estimates from \Cref{semilogalignbounds} into the strategy of Benoist and Quint. We estimate:
    \begin{align*}
        \int \exp\left( \left(\log\frac{1}{|f(x)|}\right)^\rho \right)\,d\nu(x) & \leq 1+ \sum_{n\geq 1} \nu(\{ \log |f(x)| < -(\log n)^{\rho^{-1}}\})\\
        &\leq 1 + \sum_{n\geq 1} \exp\left(-C\left\lfloor\frac{(\log n)^{\rho^{-1}}}{\ve}\right\rfloor^{.99\delta}\right)\\
        &\leq \sum_{n\geq 1} \exp\left(-C(2\ve)^{-.99\delta} (\log n)^{.99\delta\rho^{-1}}\right)+ O(1)
    \end{align*}
    where in the last step we have used that for all but finitely many $n$,
    \[ \left\lfloor \frac{(\log n)^{\rho^{-1}}}{\ve}\right\rfloor \geq \left(\frac{(\log n)^{\rho^{-1}}}{2\ve}\right)\]
    The sum in the last line converges if $.99\delta\rho^{-1} > 1$, so that we have shown the necessary bound.
\end{proof}

While it is not quite as immediate due to subtleties with the semi-logarithm functions, this essentially proves \Cref{regularitythmintro}:
\begin{proof}[Proof of \Cref{regularitythmintro}]
    Fix $\rho < \rho' < .99\delta$. It suffices to show that
    \[ \int_{\P(\KK^d)} \exp\left(\left(\log \frac{1}{|f(x)|^2}\right)^\rho\right)< \infty. \] Clearly for any $\ve_0> 0$,
    \[ \int_{f(x) \geq \ve_0} \exp\left(\left(\log \frac{1}{|f(x)|^2}\right)^\rho\right) < \infty. \]
    At the same time, there is $\ve_0 > 0$ such that for $|f(x)| < \ve_0$, we have:
    \begin{equation*}
        \exp\left(\left(\log \frac{1}{|f(x)|^2}\right)^\rho\right) 
        \leq \exp\left(\left(\log \frac{1}{|f(x)|}\right)^{\rho'}\right)
    \end{equation*}
    But integrability of the right hand side term is precisely what was shown in \Cref{semilogfunctionalreg}.
\end{proof}

\begin{remm}
    We have shown the distribution weak-H\"older, with $\rho < .99\delta$. We can in fact take $\rho < \delta$ via the arguments here; recall that $.99\delta$ can be replaced with any $\delta'$ smaller than $\delta$.
\end{remm}

\subsection{Proof of continuity of variance}

The proof of \Cref{varcont} requires us to prove continuity of the variance associated to a centered cocycle. Specifically, to any $\mu \in \mathcal{P}_{\log}^2$ there is a continuous function $\psi^\mu$ on $\P(\KK^d)$ such that
\begin{equation}
    \int_{g \in G} \Phi(g,x) - \psi^\mu(x) + \psi^\mu(gx)\,d\mu(g) = \lambda_1(\mu)
\end{equation}
and by work of Benoist and Quint, the standard deviation $\sigma_\mu$ appearing in (\ref{clt}) is precisely given by:
\begin{equation}\label{varianceCocycleEq}
    \sigma_\mu := \int_{\SL_d(\KK) \times \P(\KK^d)} (\Phi(g,x) - \psi^\mu(x) + \psi^\mu(gx) - \lambda_1(\mu))^2\,d\mu(g)\,d\nu(x)
\end{equation}

Clearly:
\begin{align*}
    \left|\sigma_\mu - \sigma_{\mu'}\right| &\leq \left|\int \Phi^2(g,x)\,d(\mu \times \nu - \mu' \times \nu')\right|\\ 
    &+ 2\| \psi^\mu - \lambda_1(\mu) -\psi^{\mu'}+ \lambda_1(\mu')\|_\infty\left|\int \Phi(g,x) \,d(\mu \times \nu - \mu' \times \nu')\right|\\ &+ \| \psi^\mu - \lambda_1(\mu) -\psi^{\mu'}+ \lambda_1(\mu')\|_\infty^2
\end{align*}
Control of the first term essentially follows by Kantorovich duality; if $\mu_n \rightarrow \mu$ in $W^2_{\log}$, then in particular this also implies weak convergence of $\nu_n$ to $\nu$. By compactness, $\nu_n \rightarrow \nu$ in the $W^1$ sense with respect to the standard Fubini-Study metric we have used throughout. For a natural choice of metric on $\SL_d(\KK) \times \P(\KK^d)$, i.e.:
\begin{equation}
    \tilde{d}((M_1,x_1),(M_2,x_2)) = \log^{2\star}(\|M_1-M_2\|) + d(x_1,x_2)
\end{equation}
we then get $\mu_n \times \nu_n \rightarrow \mu \times \nu$ in $W^1( \SL_d(\KK) \times \KK^d)$, where this is precisely the Wasserstein distance with respect to the metric $\tilde{d}$. $\Phi^2$ is Lipschitz with respect to this metric, and so by Kantorovich duality the first term vanishes as $\mu' \rightarrow \mu$.

Thus Theorem \ref{varcont} follows immediately from the following convergence results for the coboundaries $\psi$.
\begin{lem}
    If $\mu_n \rightarrow \mu$ in the $W_{\slog}^\delta$ topology for some $\delta$, then $\|\psi^{\mu_n} - \psi^\mu\|_\infty \rightarrow 0$.
\end{lem}
\begin{proof}
    First, we recall precisely how these $\psi^\mu$ are defined; for $\mu$ we let $\tilde{\mu}$ be its pushforward by the inverse map $M \mapsto M^{-1}$. Then $\tilde{\mu}$ is also strongly irreducible and contracting, and moreover $\mu \mapsto \tilde{\mu}$ is continuous as a map on all of our Wasserstein distances. So in particular, there is a unique stationary measure $\nu^\ast$ associated to $\tilde{\mu}$. Then the coboundary is defined as follows:
    \begin{equation}
        \psi(x) = \int_{\P((\KK^d)^\ast)} \log d(x,y) \,d\nu^\ast(y).
    \end{equation}
    That this converges and is continuous is by no means obvious; it was in a certain sense the central technical result in \cite{bqclt}. Letting $\mu_n$ be a sequence of measures converging to $\mu$ in the $W^{\delta}_{\slog}$ topology and $\nu^\ast_n$ the associated stationary measures for the inverted distribution, the lemma naturally reduces to
    \begin{equation}
        \int_{\P((\KK^d)^\ast)} \log d(x,y)\,d(\nu_n^\ast - \nu^\ast)(y) \rightarrow 0.
    \end{equation}

    Naturally, we split this integral into two pieces; for some $\ve_n > 0$ to be determined we consider
    \begin{equation}
        \int_{d(x,y) \geq \ve_n} \log d(x,y)\,d(\nu_n^\ast - \nu^\ast)(y)
    \end{equation}
    and
    \begin{equation}
        \int_{d(x,y) < \ve_n} \log d(x,y)\,d(\nu_n^\ast - \nu^\ast)(y)
    \end{equation}
    As long as $\ve_n \rightarrow 0$, the second integral vanishes in the limit, essentially as a consequence of the computations in the proof of \cite[Proposition 4.5]{bqclt}; in particular for some $c> 0 $ depending on $\lambda_1 - \lambda_2$ we have:
    
    \begin{equation} \label{secondintbd}\int_{d(x,y) < \ve_n} \log d(x,y)\,d(\nu_n^\ast - \nu^\ast)(y) \leq (p-1) 2^p c^{p-1} \sum_{ k\geq \lfloor \frac{1}{c\ve_n}\rfloor} n^{p-2} \P[d(x,y) \leq e^{-ck}]
    \end{equation}
    and the right hand side converges to zero uniformly as $\ve_n \rightarrow 0$. It suffices to find $\ve_n \rightarrow 0$ such that the first integral converges.

    Recall that on compact metric spaces, the $W^1$ topology coincides with the weak topology, and so in particular the fact that $\nu^\ast_n$ converges to $\nu^\ast$ weakly in fact implies that $W^1(\nu^\ast_n, \nu^\ast) \rightarrow 0$. By Kantorovich duality, this means that for any $f$ with Lipschitz constant $c_f$, we have $\int f d\,(\nu^\ast - \nu_n^\ast) \leq c_fW^1(\nu^\ast, \nu^\ast_n)$.

    Finally,
    \begin{align*}
        \left| \int_{d(x,y) \geq \ve_n} \log d(x,y)\,d(\nu^\ast - \nu^\ast_n) \right| &\leq \int_{\P((\KK^d)^\ast))} |\max\{\log d(x,y),\log \ve_n\}|\,d(\nu^\ast - \nu^\ast_n)\\
        &\leq \frac{1}{\ve_n} W^1(\nu^\ast, \nu^\ast_n).
    \end{align*}
    Taking $\ve_n = \sqrt{W^1(\nu^\ast, \nu^\ast_n)}$, we obtain the desired convergence.
\end{proof}
\begin{remm}
    We note that (\ref{secondintbd}) also holds uniformly for a compact subset of $W^2_{\log}(\SL_2(\KK))$ satisfying the necessary dynamical assumptions; the only obstruction in the general case $W^2_{\log}(\SL_d(\KK))$ (as opposed to the semi-logarithmic regime) is absence of a uniform lower bound on $\lambda_1(\mu) - \lambda_2(\mu)$, which is obvious for $\SL_2$ cocycles; $\lambda_1 - \lambda_2 = 2\lambda_1$, and $\lambda_1$ is a positive continuous quantity. Hence we also obtain a proof of \Cref{varcont2d}.

    We additionally remark that our result is entirely non-quantitative; we have no modulus of continuity. Quantitative versions of continuity would require, among other things, quantitative estimates on the $W^1$ distance between two inverse dual measures $\nu^\ast$ and $(\nu^\ast)'$ in terms of some distance on $\mu$, $\mu'$ giving rise to them; even in the bounded case using e.g. the $W^\infty$ distance defined in \Cref{geodesicsec} this seems difficult.
\end{remm}

\section{Applications}\label{applsec}

\subsection{Random Schr\"odinger operators}\label{schrodsec}
Recall that the Anderson model in one dimension is a random operator of the form
\[H = \Delta + V\]
acting on $\ell^2(\Z)$ where $V$ is a random potential acting by multiplication with $V_n$ i.i.d. with law $\mu$.

Associated to such operators are the so-called transfer matrices, introduced again for convenience:

\begin{equation}\label{transfermatrixeq}
    A_n^E:= \begin{cases}
        \begin{pmatrix}
        E - V_n & -1 \\ 1 & 0
    \end{pmatrix}\cdots \begin{pmatrix}
        E - V_1 & -1 \\ 1 & 0
    \end{pmatrix} &\quad \text{for n} >  0\\
    I_2 &\quad \text{for n} = 0\\
    \left(\begin{pmatrix}E - V_0 & -1 \\ 1 & 0 \end{pmatrix}\cdots\begin{pmatrix}E - V_{n-1} & -1 \\ 1 & 0 \end{pmatrix}\right)^{-1} &\quad \text{for n}< 0
    \end{cases}
\end{equation}
for any energy $E \in \R$. We also define
\begin{equation}
    A_{[x,y]}^E ;= 
        \begin{pmatrix}
        E - V_y &  -1 \\ 1 & 0
    \end{pmatrix}\cdots \begin{pmatrix}
        E - V_x & -1 \\ 1 & 0
    \end{pmatrix}
\end{equation}
for $y > x$; note that $A_{[x,y]}^E$ is distributed identically to $A_{y-x}^E$ as defined in \Cref{transfermatrixeq}. These are fundamental tools because any solution (either in $\ell^2(\Z)$ or purely in formal terms) of 
\begin{equation}\label{timeindschrod}
    H_\omega \psi = E\psi 
\end{equation}
satisfies

\begin{equation}
    \begin{pmatrix} \psi_{n+1} \\ \psi_n \end{pmatrix} = A_n^E(\omega) \begin{pmatrix}\psi_1 \\ \psi_0 \end{pmatrix}
\end{equation}

While there is no theorem to this effect in full generality, it is generally understood that if one obtains positivity of the Lyapunov exponent and large deviation estimates, both uniform in energy on compacta, then one should be able to extract almost sure Anderson localization. This has been a very successful approach, see e.g. \cite{avila2024schrodingeroperatorspotentialsgenerated,Jitomirskaya2019,Bucaj2017LocalizationFT}.

Positivity of the Lyapunov exponent in our context goes back to seminal work of Furstenberg, and continuity of the Lyapunov exponent implies uniform positivity; our contribution is the ability to get uniform large deviation estimates under more permissive moment conditions than has been previously obtained. We first prove some simple technical lemmas. We fix $E \in \mathbb{R}$, and define the map 
\begin{align*}
    f_E : \mathbb{R} &\longrightarrow \SL_2(\mathbb{R})\\
    x &\longmapsto \begin{pmatrix}
    x - E & -1\\
    1 & 0
    \end{pmatrix}
\end{align*}

\begin{prop}\label{prop:energyLogContinuous}
    Given a fixed measure $\mu\in \mathcal{P}_{\log}^p(\mathbb{R})$, the function $E \mapsto (f_E)_* \mu$ is continuous in the $p$-th logarithmic Wasserstein topology on $\mathcal{P}_{\log}^p(\SL_2(\mathbb{R}))$.
\end{prop}
\begin{proof}
    Given energies $E_1$ and  $E_2 \in \mathbb{R}$ with small difference, we consider the map
    $$F : x \longmapsto \left(\begin{pmatrix}
    x - E_2 & -1\\
    1 & 0
    \end{pmatrix}, \begin{pmatrix}
    x - E_1 & -1\\
    1 & 0
    \end{pmatrix}\right)$$
    and define the coupling $\gamma = F_*\mu$. We then have that
    \begin{align*}
        W_{\log}^p((f_{E_1})_*\mu, (f_{E_2})_*\mu) &\leq \int_{\SL_2(\mathbb{R}) \times \SL_2(\mathbb{R})} d^{p\star} (y_1,y_2) d\gamma (y_1,y_2)\\
        &= \int_{\mathbb{R}} d^{p\star} (F(x)^{(1)},F(x)^{(2)}) d\mu(x)\\
        &\leq \int_{\mathbb{R}} \log^{p\star} (|E_2 - E_1|) d\mu(x) \to 0
    \end{align*}
    as $|E_2 - E_1| \to 0$, proving the claim.
\end{proof}

By a near identical proof, one obtains also:
\begin{prop}\label{prop:energySemilogContinuous}
    Given a fixed measure $\mu \in \mathcal{P}^\delta_{\slog}$ for some $\delta \in (0,1)$, the function $E \mapsto (f_E)_* \mu$ is continuous in the semi-logarithmic Wasserstein topologies on $\mathcal{P}_{\slog}^\delta(\SL_2(\mathbb{R}))$.
\end{prop}

These results give us more or less immediately \Cref{semilogschrodldes}; one does need the following lemma, which is well known.
\begin{lem}
    If $\mu$ is a measure supported on at least two points, then for any $E$, we have that $(f_E)_\ast \mu$ is strongly irreducible and contracting.
\end{lem}

We also obtain the following analogue in the semi-logarithmic regime:
\begin{thm}
    Let $\mu$ be a measure satisfying the bound 
    \[\int \exp(\max\{0,  \log |x|\}^\delta)\,d\mu(x) < \infty \] for some $\delta \in (0,1)$, which is moreover not supported on a single point. Then for any $I \subset \R$ compact, and $\ve > 0$, there are $C,c > 0$ depending on $I, \ve$ and $\mu$ such that
    \begin{equation}
        \P\left[\left|\frac{1}{n}\log \left|\langle y, A_n^E x \rangle\right| - \lambda(E)\right| > \ve\right] \leq C\exp(-cn^{.99\delta})
    \end{equation}
\end{thm}

Via some modifications of existing arguments, these large deviation estimates for transfer matrices will yield exponential decay of what is called the Green's function, and also enable the Wegner estimate, the two new ingredients enabling the proof of \Cref{semiloglocal} in \cite{hurtado2025localization}.

\subsection{Random geodesics on hyperbolic surfaces}\label{geodesicsec}

We shall now apply our stability results to the problem of counting random geodesics on hyperbolic surfaces. Proving \Cref{HyperbolicCountingCLT} amounts to showing that the mean and variance of the CLT attached to the probability measure $\mu_{g}$ on $\mathrm{PSL}_2(\R)$ associated to the metric $g$ on $\Sigma$ vary continuously in an appropriate Wasserstein topology --- since such a measure is still contracting and strongly irreducible, the theorem will follow. We first describe the Teichm\"{u}ller space of hyperbolic structures and its equivalent reformulation in terms of the holonomy representation. See \cite[Chapters 10, 11]{FarbMargalit} for a more complete introduction.

\begin{deffo}
    The \textit{Teichm\"{u}ller space} $\mathcal{T}(\Sigma)$ is the space of marked hyperbolic surfaces $[X,\varphi : \Sigma \to X]$ up to homotopy (i.e. $\varphi_1 \sim \varphi_2$ when $\varphi_2 \circ \varphi_1^{-1} \simeq \mathrm{id}$).
\end{deffo}

We will frequently abuse notation and refer to a point in Teichmuller space using the hyperbolic surface $X$ or its associated metric $g$ on $\Sigma$. 

The Teichm\"{u}ller space is topologized via Fenchel--Nielsen coordinates, through which one can see it is homeomorphic to a ball of complex dimension $3g-3 +n$, where $g$ is the genus of $\Sigma$ and $n$ is the number of punctures. It can also equivalently be built as the space of discrete, faithful, type preserving representations up to conjugacy $\mathrm{DF}(\pi_1(\Sigma), \mathrm{PSL}_2(\R))/\mathrm{PSL}_2(\R)$ under the compact-open topology --- these two topologies on $\mathcal{T}(\Sigma)$ are equivalent.

Let $\mathcal{S}$ denote the set of standard symmetric generating set of the surface group $\pi_1(\Sigma)$, define $\delta_\mathcal{S}$ to be the uniform distribution on $\mathcal{S}$, and let $\delta_\mathcal{S}^{*n}$ denote the $n$-th convolution power of $\delta_\mathcal{S}$. Essentially all that is necessary is continuity of $\delta_{\mathcal{S}}$ as a function of $g$ in an appropriate Wasserstein topology; $W^\delta_{\slog}(\mathrm{SL}_2(\R))$ would suffice. However, the dependence of $\delta_{\mathcal{S}}$ on $g \in \mathcal{T}(\Sigma)$ is in fact much more regular. There is a topology finer than all the Wasserstein topologies thus far introduced. For $\mu, \mu'$ compactly supported measures on $\mathrm{SL}_d(\R)$, we define:
\begin{equation}\label{winf}
    W^\infty(\mu,\mu') := \inf_{\eta} \esssup_{(A,B)\sim \eta} \|A-B\|
\end{equation}
where $\eta$ ranges over all couplings of $\mu$ and $\mu'$; this is known as the $\infty$-Wasserstein distance.

Naturally, we let $\mathcal{P}^\infty(\mathrm{SL}_d(\R))$ denote the set of compactly supported measures, equipped with the topology induced by this metric. We can now prove the following lemma:

\begin{lem}\label{hyperbolicMeasureContinuity}
    Let $X$ and $Y$ be a two marked hyperbolic surfaces with associated representations $\rho_X$ and $\rho_Y$ lying in a compact subset $K\subset \mathcal{T}(\Sigma)$. The map from $\mathcal{T}(\Sigma)$ to $\mathcal{P}^\infty (\mathrm{PSL}_2(\R))$ given by
    \begin{align*}
        g_X \longmapsto \mu_{{g_X}} = (\rho_X)_* \delta_\mathcal{S}
    \end{align*}
    is continuous.
\end{lem}
\begin{proof}
Recall that the $\infty$-Wasserstein distance is defined by
    $$W^{\infty} (\mu_{g_X},\mu_{g_Y}) = \inf_\eta \esssup_{(A,B)\in \mathrm{PSL}(2,\R)^2} 
 \|A - B\|$$
 where $A$ and $B$ are distributed according to the projections of the coupling $\eta$ onto each factor of $\mathrm{PSL}(2,\R)$. Since the map $g_X \mapsto \mu_{g_X}$ factors as the limit of a sequence of pushforwards, and convolution of the measure $\delta_\mathcal{S}$ with itself is continuous, it suffices to check that the Wasserstein distance between the pushforward measures $(\rho_X)_* \delta_\mathcal{S}$ and $(\rho_Y)_* \delta_\mathcal{S}$ is small when $\rho_X$ and $\rho_Y$ are close. 
 
 Since $\mathrm{DF}(\pi_1(\Sigma),\mathrm{PSL}(2,\R))$ is topologized by the compact-open topology, checking that we have Wasserstein continuity amounts to verifying that
 \begin{equation} \label{funnysup}
 \sup_{h \in \mathcal{S}} \|\rho_{Y} (h) - \rho_{X}(h)\| \to 0
 \end{equation}
 as $Y$ approaches $X$ in $\mathcal{T}(\Sigma)$. This immediately follows from the definition of the compact-open topology. Indeed, since we can couple $\delta_\mathcal{S}$ with itself diagonally and pushforward this measure onto $\mathrm{PSL}_2(\R)^2$ via the map $(\rho_X \times \rho_Y)$, the supremum appearing in (\ref{funnysup}) is exactly the $W^\infty$ distance between $(\rho_X)_* \delta_\mathcal{S}$ and $(\rho_Y )_* \delta_\mathcal{S}$.
\end{proof}

\begin{remm}
\begin{enumerate}
\item    The $\infty$-Wasserstein is a very strong notion of continuity; it implies continuity of the map $\mathcal{T}(\Sigma) \to \mathcal{P}^p(\mathrm{PSL}_2(\R))$ is continuous in all weaker Wasserstein topologies, and is strictly finer than another topology on the space of compactly supported measures which has been fruitfully studied, see e.g. \cite{avila2023continuitylyapunovexponentsrandom,tall2020moduli}. These topologies coincide in certain cases of interest though; both restrict to the same topology on spaces of locally constant linear cocycles over a fixed Bernoulli shift.

\item    The compact-open topology $\mathcal{T}(\Sigma)$ is metrized by the \textit{Teichm\"{u}ller distance} $d_{\mathcal{T}}$, which informally records how similar two hyperbolic surfaces are quasiconformally. One can thus ask about the modulus of continuity of the map $g \mapsto \mu_g$ in \Cref{hyperbolicMeasureContinuity}. By the results in \cite{duarte2020large}, the map $\mu_g \mapsto L_g$ is locally H\"older continuous; in particular any quantitative version of \Cref{hyperbolicMeasureContinuity} immediately gives a (local) modulus of continuity for $g \mapsto L_g$. We are not aware of any such result.
\end{enumerate}
\end{remm}

\begin{proof}[Proof of \Cref{HyperbolicCountingCLT}]
    The pushforward measures $\mu_{g_X}$ are strongly irreducible and contracting, and thus \Cref{lyapcont} and \Cref{varcont} apply.
\end{proof}
Finally, we have stability of large deviation estimates as well, although there is an important technical caveat to keep in mind. The relationship between the norm of the representation of $\gamma$ and the length of the geodesic representative only holds if its representative is hyperbolic. In the special case where $\Sigma$ is closed and hence all elements are hyperbolic we have the following:
\begin{thm}
    Let $\Sigma$ be a closed surface admitting a hyperbolic metric, and $\mathcal{T}(\Sigma)$ its Teichm\"uller space. For any $K \subset \mathcal{T}(\Sigma)$ compact, and $\ve > 0$, there is $c = c(K,\ve)>0$ such that 
    \[ \mu_n\left(\gamma\,:\, |\ell_g(\gamma) - nL_g| > n\ve\right)\leq e^{-cn}.\]
\end{thm}
If $\Sigma$ is of finite type, we still have bounds which are stable in a sense; by work of Benoist and Quint \cite{benoist2016random}, we have $\mu_n\left(\gamma\,:\, \rho_g(\gamma) \text{ is not hyperbolic} \right) \rightarrow 0$, and the rate of convergence is independent of the choice of hyperbolic metric $g \in \mathcal{T}(\Sigma)$. Moreover, unipotent elements of the fundamental group do not admit unique geodesic representatives, so the length $\ell_g$ isn't well defined on the corresponding cusps. Hence in general we get the following:
\begin{thm}\label{FiniteTypeHyperbolic}
    Let $\Sigma$ be a surface of finite type admitting a hyperbolic metric, and $\mathcal{T}(\Sigma)$ its Teichm\"uller space. Then we have
    \[\mu_n\left(\gamma \text{ hyperbolic and }\frac{\ell_g (\gamma)  - n L_g}{\sigma_g\sqrt{n}}\in [a,b]\right) \to \frac{1}{\sqrt{2\pi}}\int_a^b e^{-x^2/2} dx.\]
    Moreover, for any $K \subset \mathcal{T}(\Sigma)$ compact, and $\ve > 0$, there is $c = c(K,\ve)>0$ such that 
    \[ \mu_n\left(\gamma\,:\, |\ell_g(\gamma) - nL_g| > n\ve\right)\leq e^{-cn} + \mu_n\left(\gamma\,:\, \rho_g(\gamma) \text{ is not hyperbolic} \right).\]
\end{thm}

\bibliography{uniflyap}
\bibliographystyle{amsalpha}
\Addresses
	
\end{document}